\documentclass[smallextended,envcountsect,]{svjour3}
\usepackage{amsmath,amssymb, amsfonts}
\usepackage{algorithm}
\usepackage{float}
\usepackage{subcaption}
 \usepackage{relsize}
\usepackage{graphicx}
\usepackage{booktabs} 
\usepackage[table]{xcolor}
\usepackage{bm}
\usepackage{setspace}
\usepackage{enumitem}
\usepackage{graphicx}
\usepackage{algorithm}
\usepackage{eufrak}
\usepackage{yfonts}
\usepackage{algpseudocode}
\usepackage{cite}
\usepackage{geometry}
\usepackage{mathtools}

\renewcommand{\vec}[1] {\ensuremath{\boldsymbol{#1}}}
 \def\r{{\mathbb R}}

\def\x{{\vec x}}

\def\u{{\vec u}}
\def\v{{\vec v}}
\def\f{{\vec f}}
\def\rv{{\vec r}}
\def\g{{\vec g}}

\def\s{{\vec s}}

\def\z{{\vec z}}
\def\e{{\vec e}}
\def\y{{\vec y}}
\def\t{{\vec{\theta}}}

\def\w{{\vec w}}

\def\gbar{\overline{\g}}
\def\gstar{\g^\star}

\def\xbar{\overline{\x}}
\def\ubar{\overline{\u}}
\def\Tbar{\overline{T}}
\def\xo{\x_o}
\def\xc{\x_c}
\def\fb{{\vec{\Delta u}}}
\def\etabf{\vec{\eta}}
\def\upsilonbf{\vec{\Upsilon}}

\begin{document}

\title{Sensitivity driven experimental design to facilitate control of dynamical systems}

\author{Joseph Hart, Bart van Bloemen Waanders, Lisa Hood, Julie Parish}
\institute{Joseph Hart \at
  Optimization and Uncertainty Quantification 
   joshart@sandia.gov\\
Bart van Bloemen Waanders \at
  Optimization and Uncertainty Quantification 
  bartv@sandia.gov 
  (corresponding author) \\
  Lisa Hood \at
 Navigation, Guidance \& Ctrl II 
  lghood@sandia.gov\\
    Julie Parish \at
    Autonomy for Hypersonics
   jparish@sandia.gov
    \and
  Sandia National Laboratories\\
  P.O. Box 5800\\
  Albuquerque, NM 87123
}

\date{}

\maketitle

\abstract{Control of nonlinear dynamical systems is a complex and multifaceted process. Essential elements of many engineering systems include high fidelity physics-based modeling, offline trajectory planning, feedback control design, and data acquisition strategies to reduce uncertainties. This article proposes an optimization centric perspective which couples these elements in a cohesive framework. We introduce a novel use of hyper-differential sensitivity analysis to understand the sensitivity of feedback controllers to parametric uncertainty in physics-based models used for trajectory planning. These sensitivities provide a foundation to define an optimal experimental design which seeks to acquire data most relevant in reducing demand on the feedback controller. Our proposed framework is illustrated on the Zermelo navigation problem and a hypersonic trajectory control problem using data from NASA's X-43 hypersonic flight tests.
}

\section{Introduction} \label{sec:intro}

Optimal control has been successfully applied to a range of engineering applications from managing chemical plants to trajectory planning of aerospace vehicles. Even though solution techniques have reached impressive levels of maturity, achieving stable and robust performance while mitigating uncertainties remains a key challenge.  For instance, trajectory planning for hypersonic vehicles utilizes primal-dual interior point solvers applied to dynamics that are discretized with time adaptive \cite{hp_adap_psopt_rao} pseudospectral methods  \cite{fahroo_ross_ps_2002}. Yet significant sources of uncertainties, including environmental conditions and aerodynamic properties, result in deviations from the open loop solution which are too great for feedback controllers to overcome. In hypersonic trajectory planning, approximating aerodynamics using wind-tunnel experiments, flight data, or high-fidelity numerical models is expensive and challenging. Consequently the goal of selecting the most informative data acquisition strategies limited to a sampling budget is an ongoing research topic.

This paper introduces post-optimality sensitivities as a metric to help determine the most effective sampling strategy.  The goal is to improve the information content of the underlying models so that a feedback controller can better track the open-loop solution.  The question of how, where, and when to sample is central and because sampling real data is expensive, a judicious approach must be implemented.  We introduce a workflow that captures high-fidelity simulation through an approximation to serve as a constraint in an optimization formulation.  Through the Implicit Function Theorem, post-optimality sensitivities can identify the most influential sources of uncertainties relevant to optimal solutions. Hyper-Differential Sensitivity Analysis (HDSA) utilizes advanced numerical linear algebra for efficient computational of these sensitivities.  The algorithmic generality of HDSA has been demonstated on a range of applications problems
\cite{HDSA,rand_gsvd_hdsa,Sunseri_2020}.  In the context of trajectory planning for hypersonic vehicles, a higher order numerical model is sampled for aerodynamic coefficients in the open-loop formulation to determine solutions that avoid saturating feedback controllers.

We appeal to Optimal Experimental Design (OED) theory to provide an algorithmic foundation for the collection of experiments, data, or measurements which best inform the decision making process. For dynamical systems this could consist of providing accurate predictions, understanding certain phenomenon, or arriving at an optimal design
\cite{Atkinson92,Pazman86,Pukelsheim93,ChalonerVerdinelli95,Ucinski05,HaberHoreshTenorio_08,HoreshHaberTenorio_10,AlexanderianGloorGhattas16}. The objective function in an OED formulation quantifies a notion of uncertainty to be minimized. Instead of the traditional alphabetic optimality criteria \cite{Atkinson92}, we introduce an HDSA criteria. In traditional inverse problem formulations, the OED goal is to find the best sensor points to reduce uncertainty in the inverse solution. The focus of this paper will be on finding sampling strategies in the service of optimal control solutions.

\begin{figure}[h]
\centering
	\includegraphics[width=0.8\textwidth]{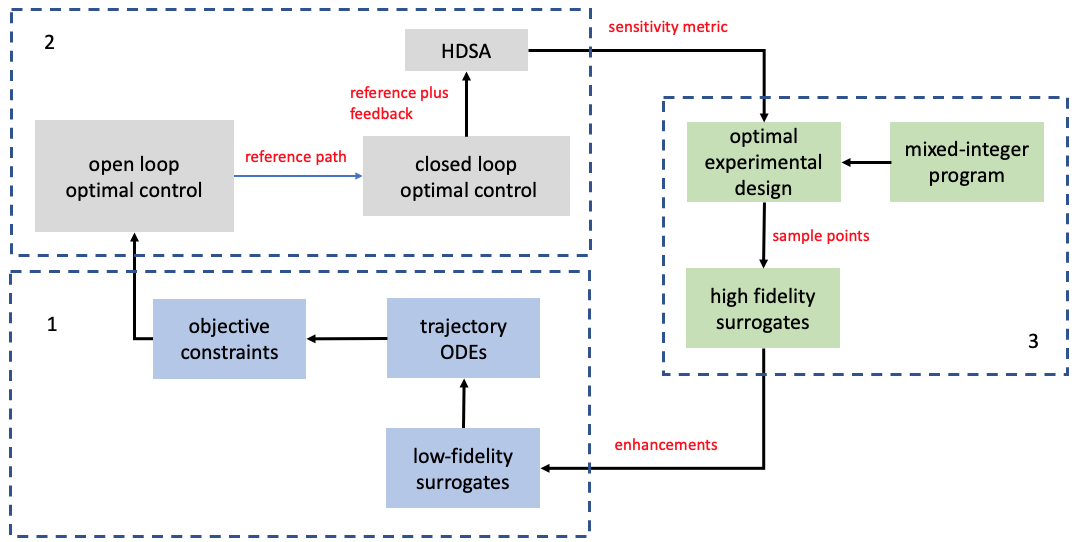}
  \caption{Analysis workflow consisting of three sequential phases: 1) trajectory ODEs supported by surrogates that feeds into an objective function, 2) open loop optimal control formulation that feeds into a closed loop optimal control problem which is instrumented to calculated HDSA, and 3) OED to determine the sampling strategy from a high fidelity surrogate.}
  \label{fig:workflow}
\end{figure}

Our algorithmic approach is organized in three phases, each consisting of supporting elements, depicted in Figure~\ref{fig:workflow}. The first phase concerns modeling. This feeds into the second phase where optimization is used to control the system. Given a model and control strategy, the third phase uses OED to direct data acquisition which seeks to reduce the modeling uncertainties most relevant to the control solution. The first phase encapsulates the dynamics through a set of ordinary differential equations (ODEs). These ODEs often contain parameters which are estimated via a table-lookup or polynomial-based surrogates based on experiments or expensive high-fidelity computations. The users specification of an objective and constraints completes the necessary modeling to pose the optimization problems in phase 2. Within this phase, there are two distinct optimization problems, open loop optimal control to generates a reference trajectory and closed loop optimal control to equip the system with feedback to mitigate the effects of uncertainty. However, closed loop control is challenging in the face of significant uncertainty so we introduce HDSA on the closed loop control problem to prioritize which uncertainties place the greatest strain on the feedback controller. This sensitivity information flows into phase 3 where we pose an OED problem which is solved using mixed-integer programming. The resulting design informs data acquisition which is incorporated into the phase 1 modeling to complete the iterative workflow.

Our contributions consist of 1) the introduction of hyper-differential sensitivities to identify sources of uncertainty that place the greatest demand on a feedback controller; 2) defining an optimal experimental design criteria which leverages the hyper-differential sensitivities to direct data acquisition most relevant to reducing demand on the feedback controller; 3) demonstrating our proposed framework to find an optimal trajectory path using dynamics that are reprensentative of a real hypersonic vehicle.  The remainder of the article is organized as follows. Section~\ref{sec:formulation} introduces a general problem formulation for nonlinear control motivated by our hypersonic flight exemplar. We present an interpretation of Hyper-differential sensitivity analysis in Section~\ref{sec:hdsa} to facilitate a connection between trajectory planning and feedback control. These sensitivities are used to develop a novel optimal experimental design in Section~\ref{sec:model_enhancement}. We illustrate these ideas in Section~\ref{sec:numerics} on the Zermelo navigation problem, an insightful prototype, and an application to the NASA X-43 hypersonic vehicle \cite{Harsha_2011}. A concluding perspective is given in Section~\ref{sec:conclusion}.

\section{Problem formulation} \label{sec:formulation}
We consider control of a nonlinear dynamical system governed by ordinary differential equations. Letting $t$ denote time, $T>0$ denote a final time, $\x: [0,T] \to \r^n$ be time dependent state variables ($\dot{\x}$ denotes their time derivative), and $\u:[0,T] \to \r^m$ denote controllers, we consider
\begin{align}\label{eq:OL_OCP}
 \tag{OL($\g$)} 
& \min_{\x,\u,T}  \int_0^{T} C_{run}(\x(t),\u(t),t)dt + C_{final}(\x(T),\u(T),T) \\
& s.t. \nonumber \\
& \begin{dcases} \label{eq:dynamics}
 \dot{\x}(t) = \f(t,\x,\u,\g(\x,\u)) \qquad t \in (0,T) \\
  \x(0) = \x_0
 \end{dcases}\\
 & \begin{dcases}
\etabf(t,\x,\u,\g(\x,\u)) \le 0 \qquad t \in (0,T)  \\
\upsilonbf(\x(T),\u(T))=0
 \end{dcases}
\end{align}
where $C_{run}$ and $C_{final}$ denote running and final time objective functions, $\etabf \in \r^r$ denotes a set of $r$ inequality constraints and $\upsilonbf\in \r^s$ denotes a set of $s$ terminal equality constraints. The constraints $\etabf$ and $\upsilonbf$ are frequently a critical part of the problem formulation where $\etabf$ ensures that the system operates in a safe regime (bounding temperature, pressure, ect) while $\upsilonbf$ enforces that the controller drive the system to the desired final state. The solution of~\eqref{eq:OL_OCP} is frequently call the open loop control strategy.

We assume that the state dynamics 
$$\f(t,\x,\u,\g(\x,\u)) \in \r^n$$
are of high fidelity but the overall model fidelity is limited by errors in
$$\g(\x,\u) \in \r^\ell $$ 
which models external forces, system properties, and/or the environment. In our motivating hypersonics application, $\f$ corresponds to the equation of motion for a vehicle and $\g$ models aerodynamic forces acting on it. 

Throughout the article we will assume that there exists a ``true" model $\gstar$ which the system experiences in practice. Rather than knowing $\gstar$, we have access to an approximation $\gbar \approx \gstar$. Solving~\eqref{eq:OL_OCP} with $\g=\gbar$ generates a reference solution which we denote as $(\xbar(t),\ubar(t),\Tbar)$. Running the system with controller $\u=\ubar$ is called open loop control and will yield a state trajectory $\xo$ which satisfies the dynamics~\eqref{eq:dynamics} with $\u=\ubar$ and $\g=\gstar$. For many systems, $\xo$ is far from $\xbar$ due to the error $\gstar-\gbar$ and hence open loop control is insufficient. To combat this challenge, closed loop control seeks to use feedback of the system state in order to adjust the control strategy. In other words, it defines a feedback law $\fb(\x)$, mapping the state to a control response, to augment $\ubar$. This yields a closed loop trajectory $\xc$ which satisfies the dynamics~\eqref{eq:dynamics} with $\u=\ubar+\fb(\x)$ and $\g=\gstar$. Table~\ref{tab:notations} summarizes this notation for the readers convenience. 

\begin{table}[!ht]
\centering
\begin{tabular}{c|c|c}
State & Control & Model  \\ \hline
$\xbar$ & $\ubar$ & $\gbar$ \\
$\xo$ & $\ubar$ & $\gstar$ \\
$\xc$ & $\ubar+\fb(\x)$ & $\gstar$ \\	    
\end{tabular}
\caption{Summary of notation for the state, control, and model. In each row, the left column denotes the solution of the ODE system~\eqref{eq:dynamics} with final time $T=\Tbar$, control $\u$ given by the center column, and model $\g$ given by the right column.}
\label{tab:notations}
\end{table}

Determining a feedback strategy $\fb(\x)$ such that $\xc \approx \xbar$ is challenging for highly nonlinear systems. Prevalent challenges include physical constraints of system which limit the computational power available for feedback, fast time scales which mandate rapid computation, and strong nonlinearities which require time-intensive computation to resolve. To combat this challenge, we consider evaluating $\gstar(\x,\u)$ at a small number of state-control inputs which must be selected judiciously in order to produce a reference solution which may be more easily tracked by a feedback controller. We propose to use hyper-differential sensitivity analysis coupled with a novel optimal experimental design formulation to direct data acquisition. This guides offline computation to enable better trajectory generation and feedback control.

\section{Hyper-differential sensitivity analysis} \label{sec:hdsa}
To determine which sources of error in $\g$ place the greatest strain on the feedback controller $\fb(\x)$, let $\e=\x-\overline{\x}$ denote the deviation of the state $\x$ from the reference trajectory $\xbar$ and consider the reference tracking optimization problem
\begin{align}\label{eq:ref_tracking}
 \tag{RT($\g$)} 
& \min_{\e,\v} \frac{1}{2} \int_0^{\Tbar} \vert \vert \e(t)\vert \vert_2^2 dt + \frac{\alpha}{2} \int_0^{\Tbar} \vert \vert \v(t) \vert \vert_2^2 dt \\
& s.t. \nonumber \\
& \begin{dcases} \label{eq:ref_track_ode}
 \dot{\e}(t) = \f(t,\xbar+\e,\ubar+ \v ,\g(\xbar+\e,\ubar+\v)) -  \f(t,\xbar,\ubar ,\g(\xbar,\ubar))  \qquad t \in (0,\Tbar) \\
 \dot{\e}(0) = \vec{0}
 \end{dcases}
\end{align}
where $\v:[0,\Tbar] \to \r^m$ augments the controller $\ubar$ and $\alpha \ge 0$. The ODE system~\eqref{eq:ref_track_ode} follows from differentiating $\e=\x-\overline{\x}$ with respect to $t$ and substituting the dynamics~\eqref{eq:dynamics}.

For a given $\g$, the solution of~\eqref{eq:ref_tracking} gives the control $\v$, augmenting $\ubar$, such that the state optimally tracks $\xbar$.  We interpret $\v$ as the optimal closed loop controller, which cannot be computed in practice since $\gstar$ is unknown, but is useful for analysis. Contrary to feedback controllers $\fb:\r^n \to \r^m$ which map the state $\x(t)$ to a control strategy $\fb(\x(t))$, $\v(t)$ is defined on the time interval $[0,\Tbar]$. This allows us to consider the optimal closed loop control without requiring a specific functional form for $\fb(\x)$, thus making our framework applicable for a wide range of feedback control strategies. To facilitate our analysis, we assume that $\v$ and $\g$ belong to Hilbert spaces (for instance, square integrable functions) and that $\f$ and $\g$ are twice continuously differentiable. 

Hyper-differential sensitivity analysis (HDSA) \cite{HDSA} considers the sensitivity of the solution of an optimization problem with respect to perturbations of parameters appearing in the model. In the context of reference tracking, HDSA provides insights about how perturbations of $\g$ change the optimal closed loop controller $\v$. To simplify our exposition, let $\vec{s}(t,\v,\g) \in \r^n$ denote the solution of the reference tracking ODE system~\eqref{eq:ref_track_ode}, which we assume to be unique for any given $(\v,\g)$. Then the reduced space problem
\begin{eqnarray}
\label{eq:ref_tracking_rs}
\min_{\v} J(\v;\g) : = \frac{1}{2} \int_0^{\Tbar} \vert \vert \vec{s}(t,\v(t),\g)\vert \vert_2^2 dt + \frac{\alpha}{2} \int_0^{\Tbar} \vert \vert \v(t) \vert \vert_2^2 dt 
\end{eqnarray}
admits the same control solution $\v$ as~\eqref{eq:ref_tracking}. Observe that $\overline{\v} \equiv 0$, the function mapping $[0,\Tbar]$ to 0, is the global minimizer for~\eqref{eq:ref_tracking_rs} when $\g=\gbar$ (since $\xbar$ was generated with $\g=\gbar$).

 Noting that, for any fixed $\g$, the derivative of $J$ with respect to $\v$, $J_\v$, equals zero at the minima of $J$, we apply the implicit function theorem to the equation
\begin{eqnarray*}
J_\v(\v;\g)=0.
\end{eqnarray*}
Since $J_\v(\overline{\v} ,\gbar)=0$, the implicit function theorem yields the existence of an operator 
$$\mathcal F:\mathcal N(\gbar) \to \mathcal N(\overline{\v} ),$$
 defined on neighborhoods of $\gbar$ and $\overline{\v}$, such that
\begin{eqnarray*}
J_\v(\mathcal F(\g);\g)=0,
\end{eqnarray*}
i.e. $\mathcal F$ maps models $\g$ to optimal controllers $\v$. Further, $\mathcal F$ is differentiable and its derivative is given by
\begin{eqnarray}
\label{eqn:D}
\mathcal F_\g(\gbar) = - \mathcal H^{-1} \mathcal B,
\end{eqnarray}
where $\mathcal H=J_{\v,\v}$ is the hessian of $J$ with respect to $\v$, and $B=J_{\v,\g}$ is the mixed second derivative of $J$ with respect to $\v$ and $\g$, both evaluated at $(\overline{\v},\gbar)$. 

Since $\mathcal F(\gstar)$ defines the optimal closed loop controller to track $\xbar$, the control effort (norm of the closed loop controller) required to track the reference trajectory may be approximated as
\begin{align*}
\vert \vert \mathcal F(\gstar) \vert \vert = \vert \vert \mathcal F(\gstar) - \mathcal F(\gbar) \vert \vert  \approx \vert \vert \mathcal F_\g(\gbar) (\gstar-\gbar) \vert \vert
\end{align*}
since $ \mathcal F(\gbar)=0$ and $ \mathcal F(\gstar) - \mathcal F(\gbar)=\mathcal F_\g(\gbar) (\gstar-\gbar)  + \mathcal O( \vert \vert \gstar-\gbar \vert \vert^2)$ by Taylor's theorem. This leads us to define the hyper-differential sensitivity for a perturbation in the direction $\delta \g$ as
\begin{eqnarray}
\label{eqn:sen_fun}
\mathcal S(\delta \g) = \frac{\vert \vert \mathcal F_\g(\gbar)  \delta \g \vert \vert}{\vert \vert \delta \g \vert \vert} .
\end{eqnarray}
We interpret $\mathcal S(\delta \g)$ as the change in the optimal closed loop controller $\v$, augmenting $\ubar$, when $\gbar$ is perturbed in the direction $\delta \g$.

Directions $\delta \g$ which maximize $\mathcal S$ correspond to errors in $\gbar$ will require the greatest feedback effort to overcome. Our proposed approach finds these high sensitivity directions to inform an optimal experimental design which seeks to improve $\gbar$ in these directions, thus reducing demand on the feedback controller. In the context of our hypersonic application, $\mathcal S(\delta \g)$ identifies which flight configurations (a combination of altitude, velocity, and control surface geometry) will place the greatest stress on a feedback controller seeking to track the reference flight trajectory.

\subsubsection*{Discretization}
Since $\g$ is infinite dimensional (a function), we discretize perturbations of $\overline{\g}$ and represent them with a finite dimensional vector of parameters $\t$. For simplicity of the exposition, we focus on a single component of $\overline{\g}(\x,\u)$, without loss of generality denote it by $\overline{g}_k(\x,\u)$. In practice, the procedure described below may be done component-wise for a vector-valued $\overline{\g}(\x,\u)$. Since $\overline{g}_k$ is a model which depends on $(\x,\u)$ we parameterize a perturbation of $\overline{g}_k$ in the form
\begin{eqnarray*}
p_k(\t^k;\x,\u) = 1 + \sum\limits_{i=1}^{L_k} \theta_i^k \phi_i^k(\x,\u),
\end{eqnarray*}
where $\{\phi_i^k\}_{i=1}^{L_k}$ is a suitable basis (chosen by the user) and $\t^k=(\theta_1^k,\theta_2^k,\dots,\theta_{L_k}^k)^T \in \r^{L_k}$ is a vector parameterizing the perturbation. In practice, $\phi_i^k$ may only depend on a subset of $(\x,\u)$. We retain all input arguments for generality but $\phi_i^k$ may be constant with respect to some of the input arguments. As an example, if $\overline{g}_k$ only depends on the first state and first control variables, $x_1$ and $u_1$, then we may define a two dimensional mesh over a range of $x_1$ and $u_1$ values and define $\{\phi_i^k\}_{i=1}^{L_k}$ as local basis functions on this mesh. Then each component of $\t$ corresponds to perturbing the model $g_k$ in a specific neighborhood of $(x_1,u_1)$ values. Because of the differentiability assumptions in HDSA, we require that $\phi_i^k$ is twice continuously differentiable and suggest splines or radial basis functions as a useful discretization to accommodate smoothness and locality.

In general, $\g(\x,\u) \in \r^\ell$ so we concatenate the vectors $\t^k$ for each component to define
\begin{align*}
\t =
\left( \begin{array}{c}
\t^1 \\
\t^2 \\
\vdots \\
\t^\ell \\
\end{array} \right) \in  \r^{q}
\end{align*}
where $q=L_1+L_2+\cdots+L_\ell$.  Letting $P(\t) \in \r^{\ell \times \ell}$ be a diagonal matrix whose $(k,k)$ entry is $p_k(\t^k;\x,\u)$, we observe that 
$$P(\t) \gbar(\x,\u)$$
parameterizes perturbations of $\gbar$ such that $\t=\vec{0}$ corresponds to the nominal model $\gbar$.

Upon discretization of the control problem, which we summarize in Appendix A, the sensitivity operator $\mathcal F_\g(\gbar)$ is approximated by the matrix $D \in \r^{K_m \times q}$, where $K_m$ is the dimension of the discretized control variables. If $K_m$ or $q$ is a moderate size (for instance $\mathcal O(10)$), we may construct and store $D$ explicitly. For large-scale applications where this is not possible, $D$ may be approximated using a Singular Value Decomposition. This approach is efficient for large problems which admit low rank structure, as is commonly observed in engineering applications. We refer the reader to \cite{HDSA,rand_gsvd_hdsa} for more details.

\section{Using HDSA to inform optimal experimental design} \label{sec:model_enhancement}
Our goal is to use the sensitivity information in $D$ to determine points $\{(\x_i,\u_i)\}$ where we will evaluate $\gstar(\x_i,\u_i)$. By choosing points to reduce the error in the directions of greatest sensitivity, we generate a reference trajectory which the feedback controller may track more effectively. However, acquiring such data comes at a significant computational or experimental cost which is prohibitive. In the context of our motivating hypersonics example, the components of $\g$ correspond to aerodynamic forces and are estimated by constructing surrogate models based on wind tunnel experiments or computational fluid dynamics simulations. As the aerodynamics depend on the vehicle geometry, controller positions, and Mach number, we have limited evaluations of $\g$ relative to its complex dependencies on many variables and parameters. As the components of $\t$ correspond to our estimation of $\g$ at a particular configuration (geometry and Mach number), and the systems varies over many configurations in flight, it is imperative that we collect data in the most important configurations. A poor experimental design results in unproductive compute time and experimental efforts, and ultimately a delay or failure to yield a successful controller.

The proposed experimental design does not treat the uncertainties statistically, as it common in optimal experimental design for inverse problems \cite{alen_linear_oed,alen_nonlinear_oed}. Rather, we work in a deterministic framework focusing on the error $\gstar-\gbar$ and consider which components of this error create the greatest demand on the closed loop controller. For this reason, we will use the terminology error instead of uncertainty. 

Our optimal experimental design criteria depends on three quantities, the hyper-differential sensitivity, the expected error in the model $\g$, and the informativeness of the experiments. These three quantities, represented by matrices $D$, $A$, and $B$, will be combined to define an optimality criteria. The sensitivity information in $D$ couples the dynamics, nominal model, and reference tracking objective, while the matrices $A$ and $B$, introduced below to represent error and information gain, permit the user to incorporate domain expertise when available. Decomposing the control problem, model error estimation, and information gain from experiments enables systematic integration of application specific insights with mathematical rigor in the optimality criteria.

Let $A \in \r^{q \times q}$ be a diagonal matrix whose $(i,i)$ entry is an estimation of the error in the $i^{th}$ component of $\t$. In many cases we have some a priori insight into the relative error between components of $\t$ and this information may be embedded in $A$ so that the experimental design may exploit it. For instance, there may be greater uncertainty in aerodynamics at higher Mach numbers. If nothing is known about the error $\gstar - \gbar$ a priori, we may take $A$ as the identity matrix. The resulting optimality criteria will be unaffected by scaling $A$ by a constant, so capturing the magnitude of the error is irrelevant. Rather, the information we define in $A$ is the magnitude of the error in the components of $\t$ relative to one another.

Performing an experiment (computational or physical) will give an improvement in our estimate of some components of $\t$. We represent this by defining a vector $\w = (w_1,w_2,\dots,w_d)^T \in \{0,1\}^d$ where $d$ is the number of possible experiments and the entries of $\w$ are 1 if we conduct the experiment and 0 if we do not. We define 
\begin{eqnarray*}
r_{i,j} \in [0,1] \qquad i=1,2,\dots,q, j=1,2,\dots,d
\end{eqnarray*}
 as the relative reduction in error in the $i^{th}$ component of $\t$ when the $j^{th}$ experiment is conducted. Assuming that components of $\t$ correspond to state-control pairs $(\x_i,\u_i)$, $i=1,2,\dots,d$, a useful error reduction model is
 \begin{eqnarray*}
r_{i,j} = \exp{\left( -\gamma_j^2 \vert \vert (\x_i,\u_i)-(\x_j,\u_j) \vert \vert^2 \right)} - \epsilon_j 
\end{eqnarray*}
where $\gamma_j>0$ is a correlation length parameter and $\epsilon_j \ge 0$ is the observation error which may be zero if high-fidelity experiments are performed, or may be adapted to model experiments of varying fidelities or noise levels. When no domain expertise is available, $\gamma_j$ may be specified based on the smoothness of $\gbar$; however, when appropriate the user may impose additional domain knowledge by the specification of it. 
 
To model the error reduction provided by a set of experiments, define $B(\w) \in \r^{q \times q}$ as the diagonal matrix whose $(i,i)$ entry is given by
\begin{eqnarray*}
(B(\w))_{i,i} =b_i(\w)= \left( 1 - \sum\limits_{j=1}^d w_j r_{i,j} \right)_+ 
\end{eqnarray*}
where $(x)_+=\max\{x,0\}$ is the plus function. This models our ability to inform the $i^{th}$ component of $\t$ through a combination of experiments such that $b_i(\w) \in [0,1]$ is the relative reduction in error after conducting the experiments defined by $\w$.

To the end that collecting new data reduces the demand on the feedback controller, we would like to minimize the hyper-differential sensitivities. Ideally, for any given design we should recompute the reference trajectory and sensitivities to measure the reduction in sensitivity by the design. This is computationally intractable, and rather we use the matrix $D$ which is computed once using the nominal $\gbar$. To model the decrease in sensitivity, we consider the updated sensitivity matrix
\begin{eqnarray*}
S(\w) = D A B(\w) \in \r^{K_m \times q} 
\end{eqnarray*}
which imbeds our estimation of error  $A$ and its reduction $B(\w)$ given the experiment $\w$. We interpret $S(\w)$ as an approximation of the hyper-differential sensitivities given that the experiments defined by $\w$ are conducted. 

Let $\kappa_j > 0$, $j=1,2,\dots,d$ denote the cost of experiment $j$ and $\kappa_B > 0$ denote the experimental budget. We pose an optimization problem to minimize the sensitivity matrix $S(\w)$ given the budget constraint $\kappa_B$. With an objective of minimizing the sensitivity in an average sense, and a desire to achieve computational efficiency, we adopt the Frobenius norm (squared for mathematical convenience). This yields the optimization problem
\begin{align}
\label{eq:oed_opt}
& \min_{\w \in \{0,1\}^d} \text{Tr}(S(\w)^TS(\w)) \\
& s.t. \nonumber \\
& \sum_{j=1}^d \kappa_j w_j \le \kappa_B \nonumber 
\end{align}
where $\text{Tr}$ denotes the matrix trace. A solution of~\eqref{eq:oed_opt} is called an optimal experimental design (including the possibility of multiple feasible designs yielding the same objective function value).

We observe that~\eqref{eq:oed_opt} is an integer program with a linear constraint. Inspection of the objective function reveals that we may reformulate~\eqref{eq:oed_opt} into a standard mixed-integer quadratic program. In particular, recalling that matrix multiplication commutes in the trace, $B(\w)$ and $A$ are diagonal matrices, and
$$S(w)^TS(\w) = B(\w)^T A^T D^T D A B(\w),$$
we have
\begin{eqnarray} \label{eq:oed_obj}
\text{Tr}(S(\w)^TS(\w)) = \sum\limits_{i=1}^q d_i^2 a_i^2 b_i(\w)^2
\end{eqnarray}
where $d_i$ and $a_i$ are the diagonal entries of $D^T D$ and $A$, respectively. The scalar $d_i$ corresponds to the sensitivity function~\eqref{eqn:sen_fun} evaluated at the $i^{th}$ basis function. Physically, we may interpret~\eqref{eq:oed_obj} as scaling the sensitivity $d_i$ by the uncertainty level $a_i$ and the expected information gain $b_i(\w)$. Letting $\rv_i=(r_{i,1},r_{i,2},\dots,r_{i,d})^T$ and $c_i = d_i a_i$, $i=1,2,\dots,q$,~\eqref{eq:oed_opt} is equivalent to
\begin{align}
\label{eq:oed_opt_2}
& \min_{\w \in \{0,1\}^d} \sum\limits_{i=1}^q c_i^2 \left( 1 - \w^T \rv_i \right)_+^2 \\
& s.t. \nonumber \\
& \sum_{j=1}^d \kappa_j w_j \le \kappa_B \nonumber .
\end{align}
Since the entries of $\rv_i$, $i=1,2,\dots,q$, are nonnegative we may rewrite~\eqref{eq:oed_opt_2} as
\begin{align}
\label{eq:oed_opt_3}
& \min_{\w \in \{0,1\}^d, \s \in \r^q} \sum\limits_{i=1}^q c_i^2 s_i^2 \\
& s.t. \nonumber \\
& \sum_{j=1}^d \kappa_j w_j \le \kappa_B \nonumber \\
& s_i \ge 0, \qquad i = 1,2,\dots,q \nonumber \\
& s_i \ge 1-\w^T \rv_i,  \qquad i = 1,2,\dots,q \nonumber
\end{align}
where $\s = (s_1,s_2,\dots,s_q)^T \in \r^q$ denotes slack variables which allow us to move the plus function into inequality constraints. This poses~\eqref{eq:oed_opt_3} in a standard form for a mixed integer quadratic linear program.

Notice that the coefficients $c_i$, $i=1,2,\dots,q$, may be computed a-priori using $D$ (or a low rank approximation of it) and $A$ (which is user defined). Hence~\eqref{eq:oed_opt_3} may be fully specified to an integer programming code without needing to access original control problem or any of its operators. The results in this article are generated using Gurobi \cite{gurobi} to solve~\eqref{eq:oed_opt_3}.

\section{Numerical results} \label{sec:numerics}
In this section we provide two numerical results to illustrate important properties of the proposed method and demonstrate its effectiveness. Subsection~\ref{ssec:zermelo} provides an illustrative example using the Zermelo problem. We demonstrate the method on a hypersonic trajectory control application using aerodynamics from the NASA X43 hypersonic vehicle in Subsection~\ref{ssec:X43}.

Algorithm~\ref{alg:overview} provides a summary of the proposed approach. In it, we let $\gbar$ denote our nominal model, $(\xbar,\ubar,\Tbar)$ denote the reference solution generated using $\g=\gbar$, $\tilde{\g}$ denote our improved model after acquiring new data to improve the model, and $(\tilde{\x},\tilde{\u},\tilde{T})$ denote the updated reference solution generated using $\g=\tilde{\g}$. Our goal is that the reference trajectory $\tilde{\x}$ is easier to track than $\xbar$ since the additional data we collect reduces error in the directions which place the greatest strain on the feedback controller.

\begin{algorithm} 
\caption{Computation and use of HDSA for experimental design}
\begin{algorithmic}[1]
\State Solve the open loop problem~\ref{eq:OL_OCP} with $\g=\gbar$ to generate a reference solution $(\xbar,\ubar,\Tbar)$
\State Compute the hyper-differential sensitivities for the reference tracking problem~\eqref{eq:ref_tracking}
\State Determine the optimal experimental design by solving~\eqref{eq:oed_opt_3}
\State Evaluate $\gstar$ at the design points and use this data to fit an improved model $\tilde{\g}$
\State Solve the open loop problem~\ref{eq:OL_OCP} with $\g=\tilde{\g}$ to generate a reference solution $(\xbar,\ubar,\Tbar)$
\State Design a feedback controller to track $\tilde{\x}$
\end{algorithmic}
\label{alg:overview}
\end{algorithm}

\subsection{Zermelo problem} \label{ssec:zermelo}
\begin{figure}[h]
\centering
	\includegraphics[width=0.5\textwidth]{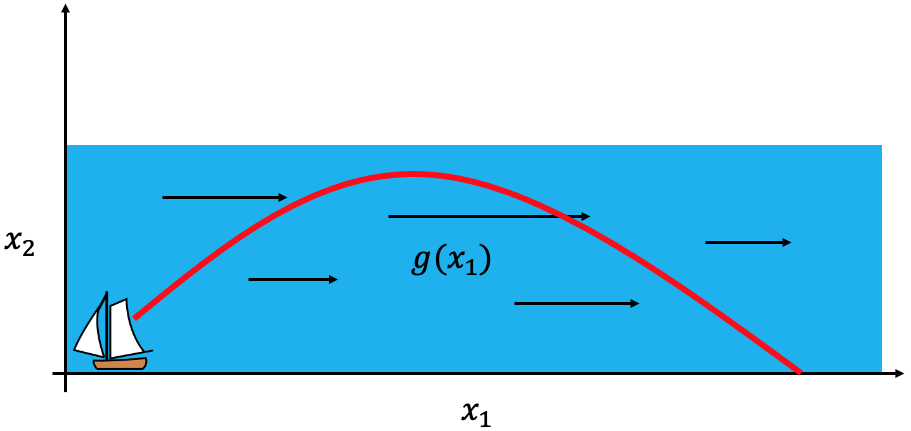}
        \includegraphics[width=0.4\textwidth]{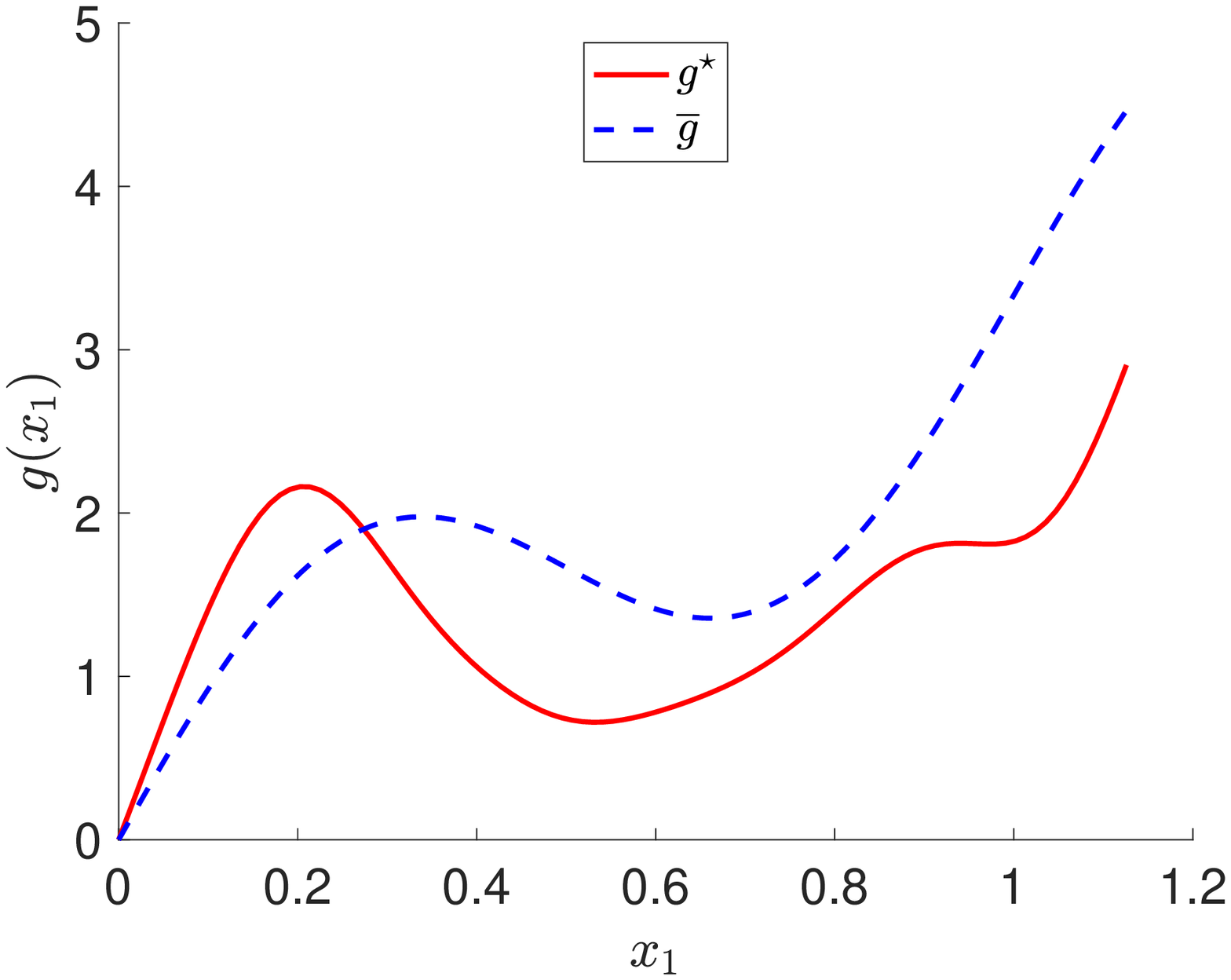}
  \caption{Left: depiction of Zermelo problem setup; right: true model $g^\star$ and initial model estimate $\overline{g}$ for the current.}
  \label{fig:velocity_fields}
\end{figure}

The Zermelo problem \cite{zermelo_example}, a classical navigation problem in the optimal control literature, considers control of a boat being driven by a current. Mathematically, it is modeled by an ODE system with states $\x=(x_1,x_2)$ corresponding to the position of a boat along a river (with the river running parallel to the $x_1$-axis). The left panel of Figure~\ref{fig:velocity_fields} depicts the problem setup. The optimal control problem under consideration is
\begin{align}\label{eq:zermelo_OC}
& \min_{\x,\u}  -x_1(1) \\
& s.t. \nonumber \\
& \begin{dcases} \label{eq:zermelo_ode}
 \dot{x}_1(t) = \cos(u(t)) + g(x_1(t))x_2(t) \qquad & t \in (0,1) \\
 \dot{x}_2(t) = \sin(u(t)) & t \in (0,1)  \\
 x_1(0) = 0 \\
 x_2(0) = 0
 \end{dcases}\\
 & \begin{dcases} \nonumber
-\pi \le u(t) \le \pi \qquad t \in (0,1) \\
x_2(1) = 0
 \end{dcases}
\end{align}
where $u:[0,1]\to \r$ is the heading of the boat which is driven by the current $g(x_1(t))x_2(t)$. This corresponds to the boat going out into the river $(x_2>0)$ and being driven downstream (in the positive $x_1$ direction) and then navigating back to the river bank $x_2=0$ at time $t=1$. The optimization objective is to maximum the range downstream.

We consider uncertainty in the current $g(x_1)$ which we seek to reduce in order to develop a robust control strategy. The right panel of Figure~\ref{fig:velocity_fields} shows the ``true" $g$, denoted as $g^\star$, and our initial estimate of it, denoted as $\overline{g}$. We assume that $g^\star$ may be evaluated at a limited number of $x_1$ values to improve $\overline{g}$. Our goal is to determine where to evaluate $g^\star$ so that we may ensure control of the system. This emulates the common application scenario where model parameters are surrogates constructed from limited evaluations of high fidelity models or controlled experiments. 

\subsubsection*{Trajectory and controller}
\begin{figure}[h]
\centering
  \includegraphics[width=0.4\textwidth]{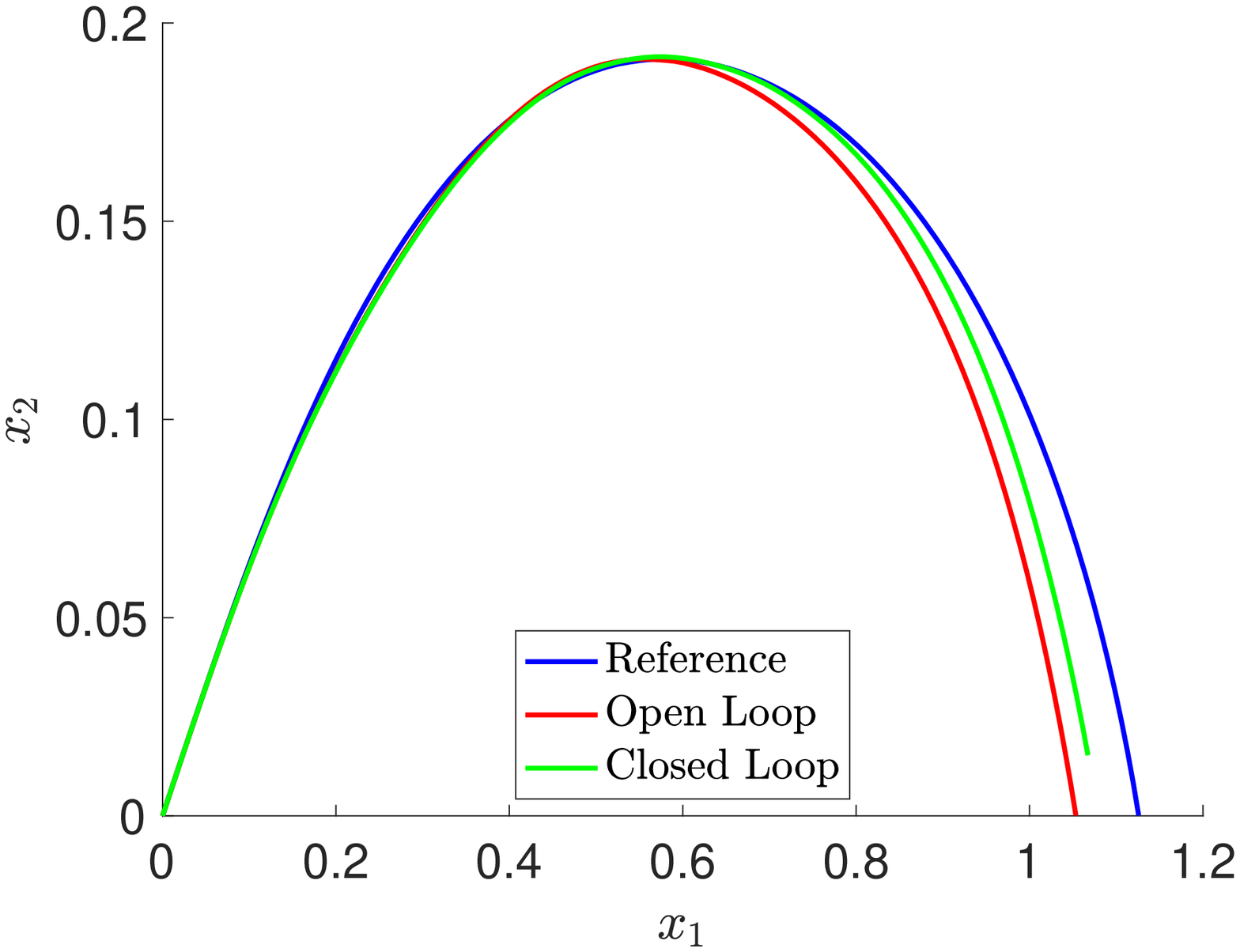}
    \includegraphics[width=0.4\textwidth]{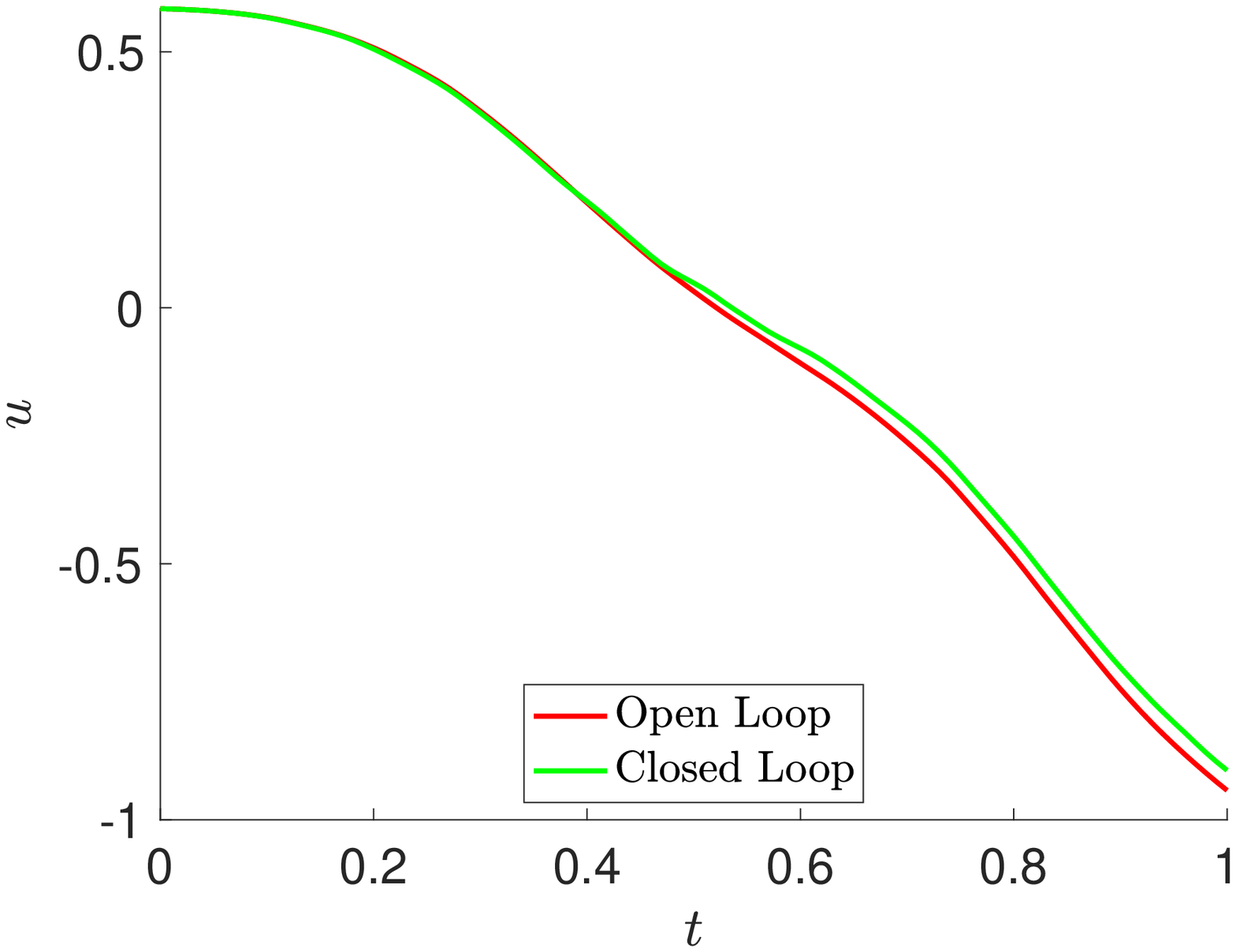}
    \caption{Trajectories in the $x_1 x_2$ phase space (left) and controllers in the time domain (right) with the reference solution $(\xbar,\overline{u})$ generated using the nominal current model $\overline{g}$. The reference trajectory $\xbar$ is the solution of~\eqref{eq:zermelo_OC} with $g=\overline{g}$ while the open and closed loop trajectories are generated using the controllers $u$ and $u+\Delta u$ with the current model $g^\star$.}
  \label{fig:trajectory_1}
\end{figure}

We solve the trajectory optimization problem~\eqref{eq:zermelo_OC} with $\overline{g}$ to generate a reference trajectory and controller $(\overline{\x},\overline{u})$. A LQR (linear quadratic regulator) feedback controller is used, where the feedback gain is computed by linearizing the dynamics around the nominal trajectory and solving a Riccati equation. The left panel of Figure~\ref{fig:trajectory_1} displays the reference trajectory $\overline{\x}$, the open loop trajectory $\xo$, and the closed loop trajectory $\xc$. The right panel of Figure~\ref{fig:trajectory_1} shows the reference (open loop) and feedback (closed loop) controllers, $\overline{u}$ and $\overline{u}+\Delta u$.

The open loop trajectory undershoots the reference trajectory as a result of $\overline{g}$ overestimating $g^\star$ for larger values of $x_1$. The feedback reduces this overshoot, but is unable to successfully track the reference trajectory. 

\subsubsection*{Hyper-differential sensitivities}
\begin{figure}[h]
\centering
        \includegraphics[width=0.4\textwidth]{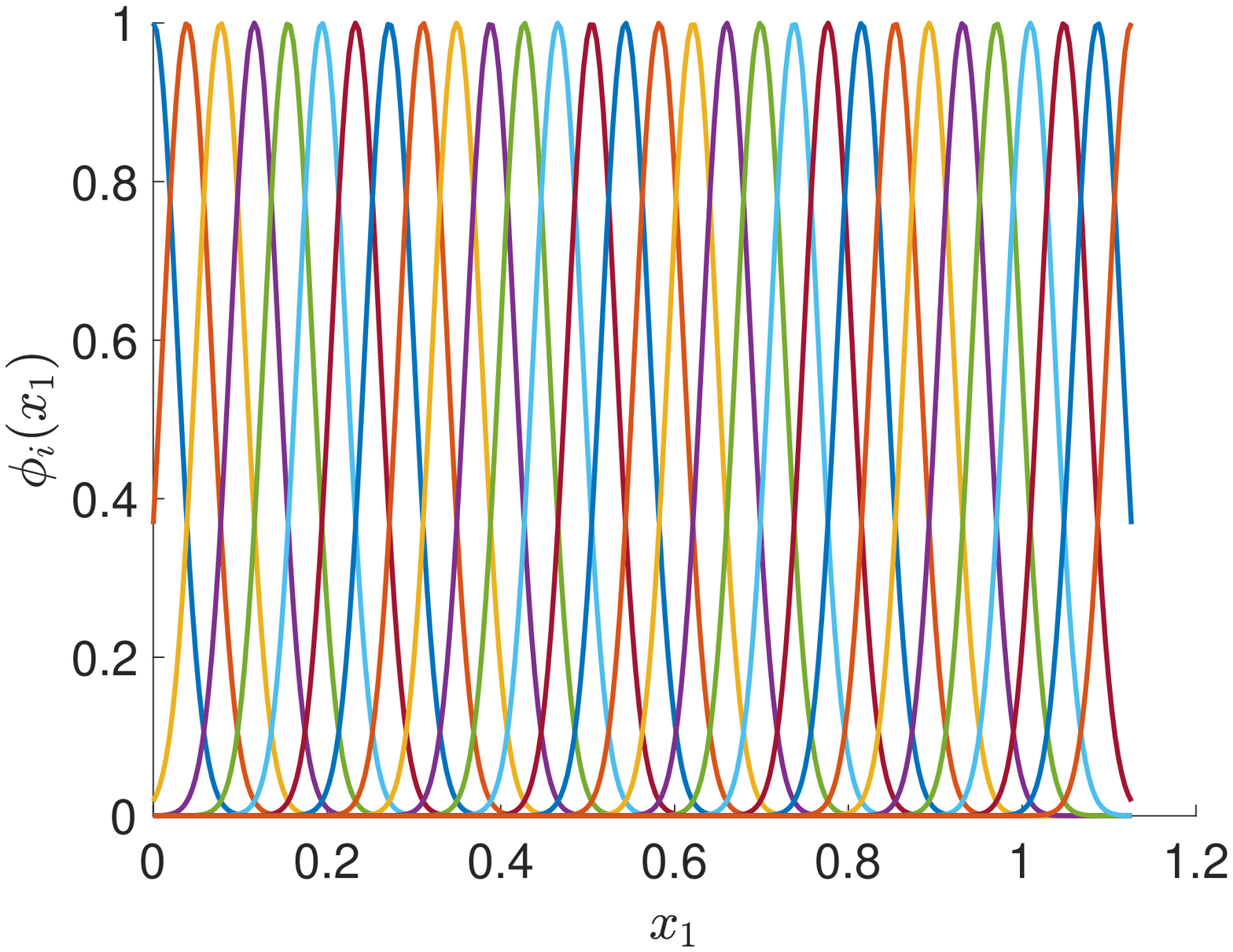}
         \includegraphics[width=0.4\textwidth]{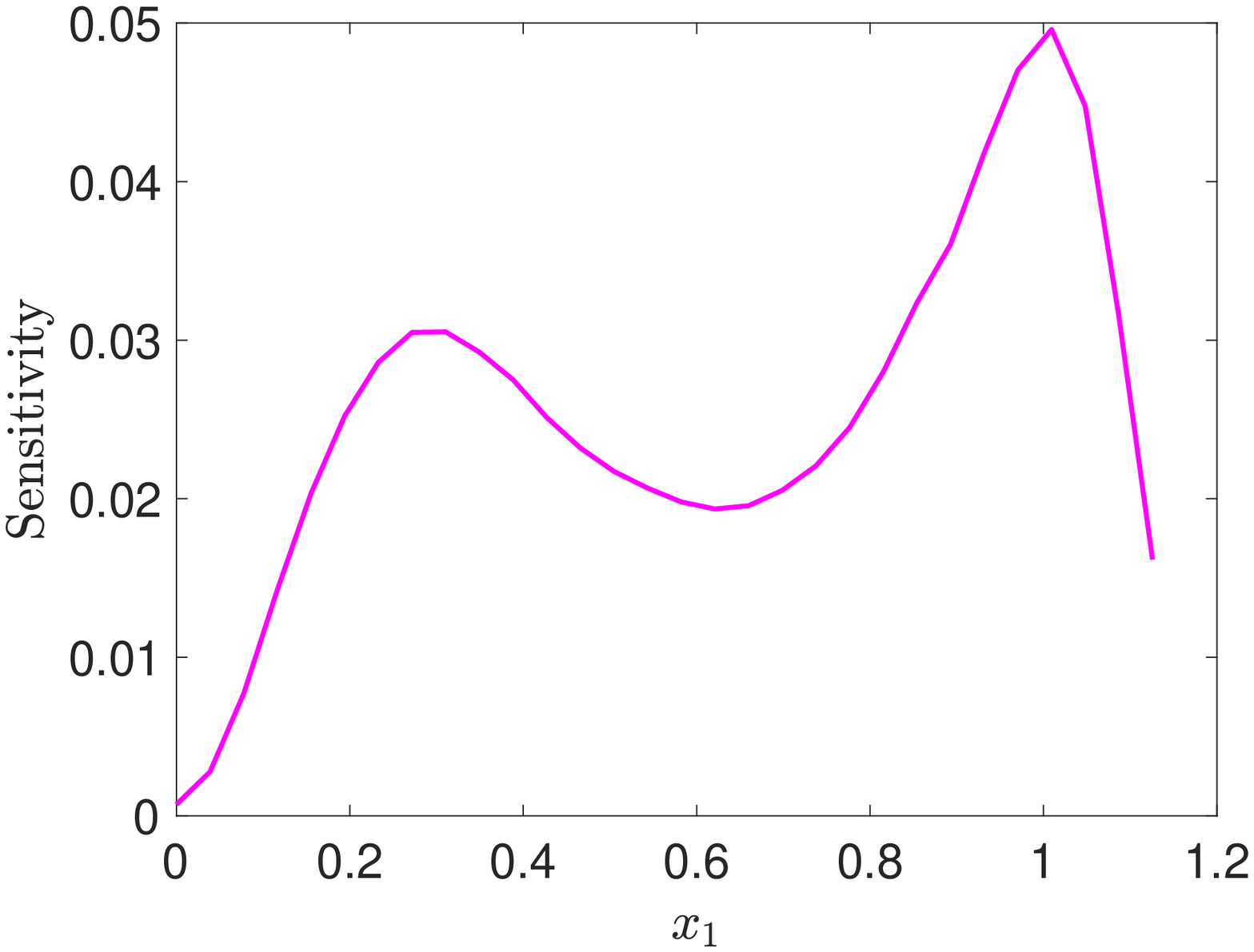}
  \caption{Left: basis functions $\{\phi_i\}_{i=1}^{30}$ used to discretize perturbations of $\overline{g}$; right: hyper-differential sensitivities for the heading control $u$ with respect to perturbations of the model $\overline{g}$.}
  \label{fig:basis_funs}
\end{figure}

We adopt radial basis functions $\phi_i(x_1)=\exp(-w(x_1-c_i)^2)$, $i=1,2,\dots,30$, which are depicted in the left panel of Figure~\ref{fig:basis_funs}, to discretize $g$. The resulting sensitivities, as a function of $x_1$, are shown in the right panel of Figure~\ref{fig:basis_funs}.

\begin{figure}[h]
\centering
         \includegraphics[width=0.4\textwidth]{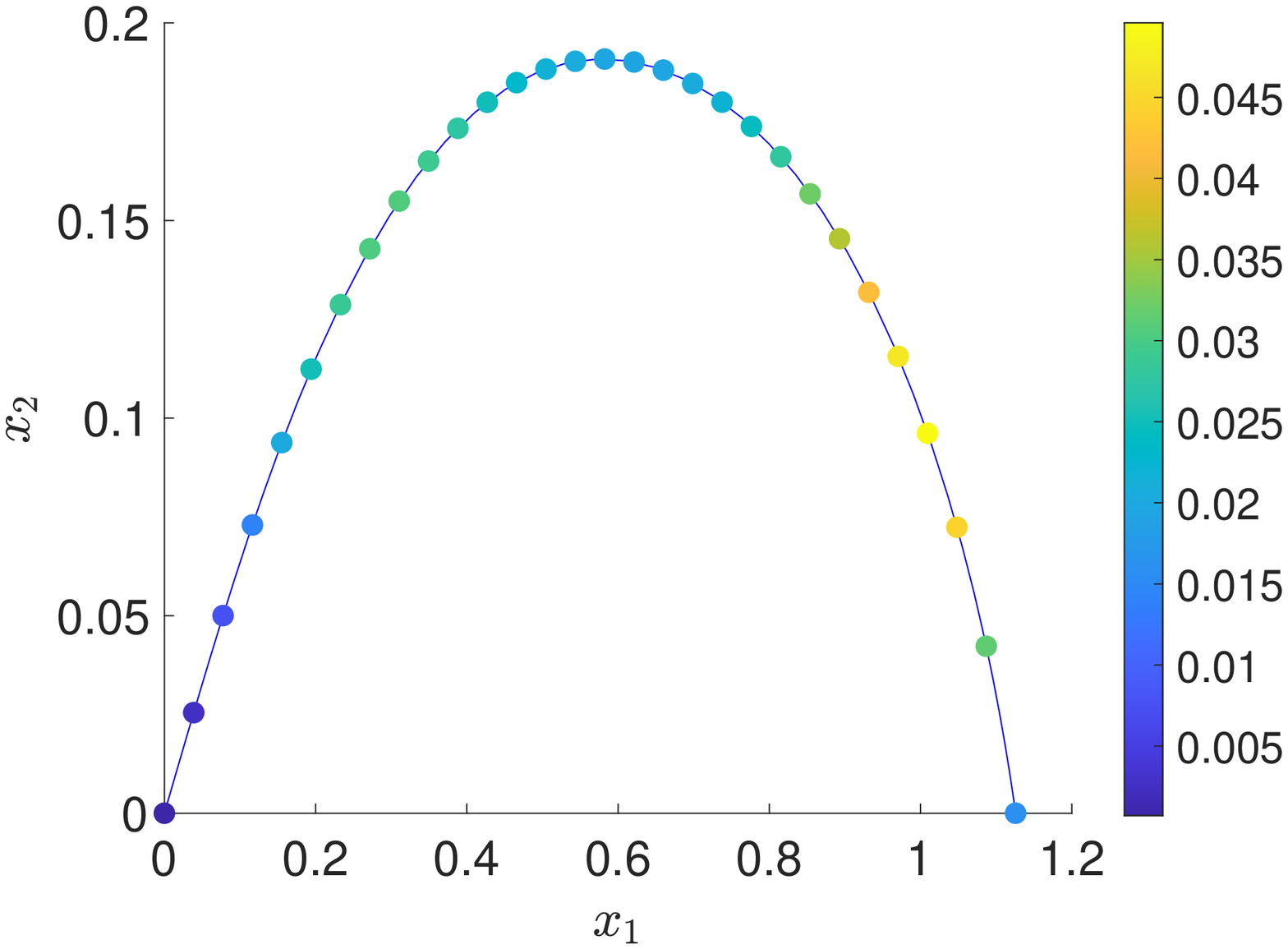}
        \includegraphics[width=0.4\textwidth]{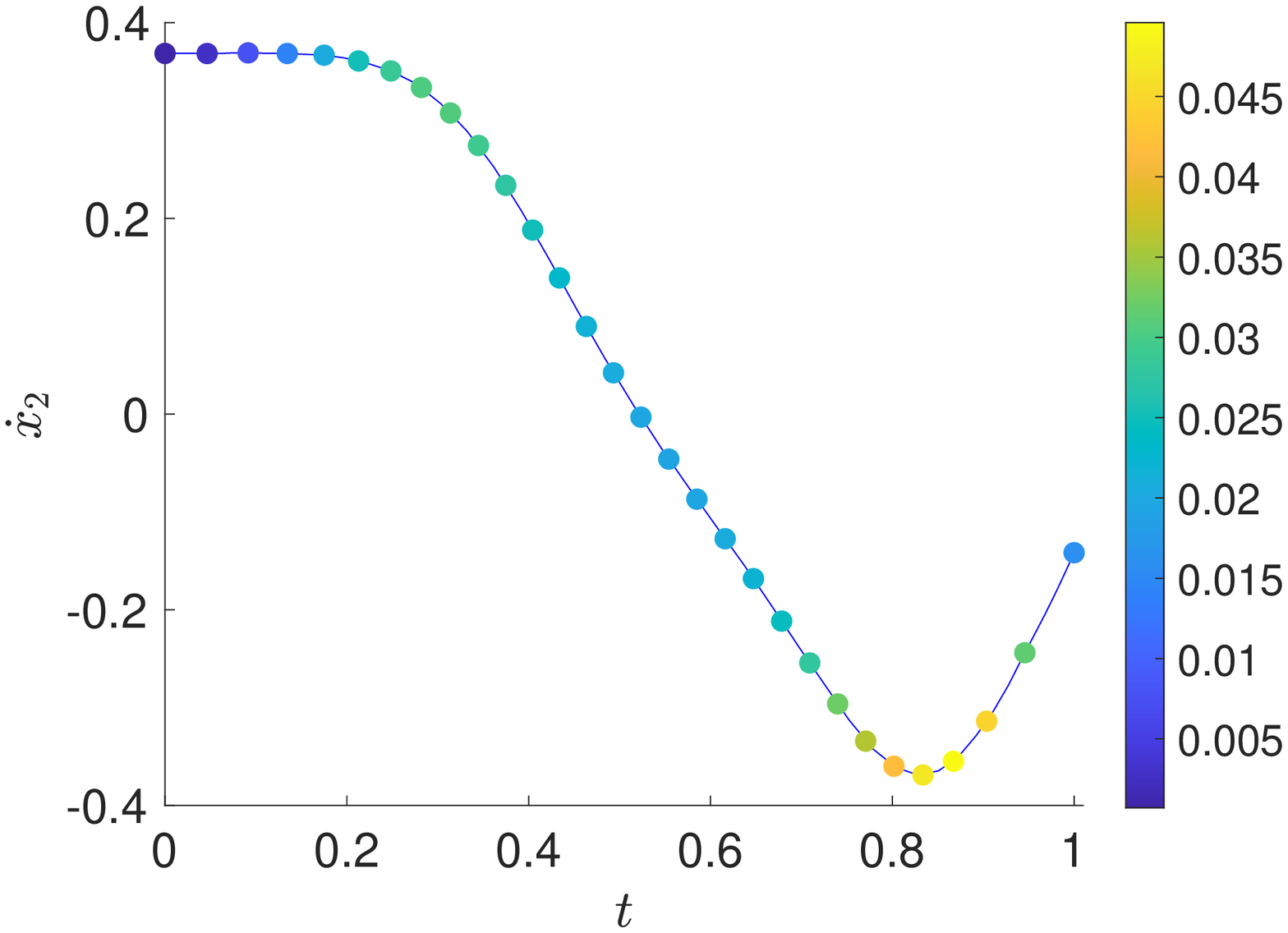}
  \caption{Left: reference trajectory $\overline{\x}$ in the $x_1x_2$ phase space; right: reference $x_2$ velocity, $\dot{\overline{x}}_2(t)$, in the time domain. Both plots have the hyper-differential sensitivities overlaid with dots; the color bar denotes the magnitude of the sensitivities.}
  \label{fig:sensitivities}
\end{figure}

To gain greater physical insight from these sensitivities, Figure~\ref{fig:sensitivities} plots $\xbar$ in the $x_1x_2$ phase space (left) and $\dot{\overline{x}}_2(t)$ in the time domain (right) with the sensitivities overlaid with colored dots (the color denoting the sensitivity magnitude). We observe that the highest sensitivity occurs when $\dot{x}_2$ is nearby its minimum. This portion of the trajectory, around $t=0.85$, is when the boat is aggressively maneuvering back to the shore as the large magnitude negative $x_2$ velocity is rapidly pushing the boat back toward the shore. The region of second greatest sensitivity is around $t=0.3$ when the boat is reducing its $x_2$ velocity so that it does not float too far out into the river (in the $x_2$ direction). We observe that these high sensitivity regions correspond to areas where $\dot{x}_2(t)$ exhibits is greatest nonlinearity as a function of $t$. This highlights the state coupling as the sensitivities correspond to a discretization of $g(x_1)$, but the high sensitivity regions correspond to nonlinearities in $\dot{x}_2$.

\subsubsection*{Experimental design}
Following the formulation of Section~\ref{sec:model_enhancement}, we seek an experimental design to determine $\kappa_B=3$ design points for which we will evaluate $g^\star$. We take a uniform mesh on $x_1$ with 30 nodes, $A$ as the identity matrix (reflecting no prior knowledge of the error $g^\star-\overline{g}$), and a error reduction model 
$$r_{i,j} = \exp(-\gamma_i^2 (x_1(w_i)-x_1(w_j))^2),$$ where $x_1(w_k)$ denotes the $k^{th}$ node in the $x_1$ mesh and $\gamma_i$ is a user defined correlation length parameter determined by scaling a nominal value by the integral of $\phi_i \vert \frac{\partial \overline{g}}{\partial x_1} \vert$ so that it is informed by the smoothness of $\overline{g}$. With an experiment cost $\kappa_j=1$ for all $j$ and budget $\kappa_B=3$, we solve~\eqref{eq:oed_opt_3} to determine an experimental design. 

\begin{figure}[h]
\centering
        \includegraphics[width=0.4\textwidth]{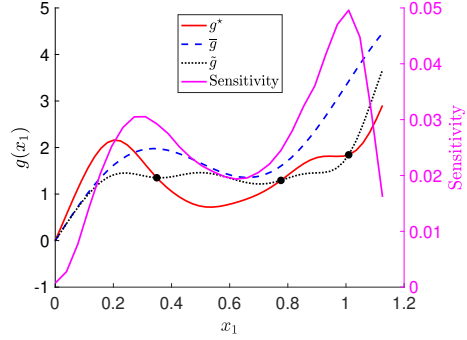}
  \caption{True, initially estimated, and improved models $g^\star$, $\overline{g}$, and $\tilde{g}$ (left axis), with experimental design points, denoted by black dots, and hyper-differential sensitivities (right axis).}
  \label{fig:improved_velocity}
\end{figure}

Using evaluations of $g^\star$ at these three points, we define the improved model $\tilde{g}=\overline{g}+\sum_{j=1}^3 \beta_j \psi_j$, where the $\psi_j$'s are radial basis functions centered at the experimental design points and the $\beta_j$'s are chosen so that $\tilde{g}$ interpolates $g^\star$ at the design points.The improved surrogate $\tilde{g}$ is displayed in Figure~\ref{fig:improved_velocity} alongside $g^\star$ and $\overline{g}$. The sensitivities from Figure~\ref{fig:basis_funs} are overlaid in Figure~\ref{fig:improved_velocity} to highlight the position of the experimental design points (indicated by black dots) corresponding to the high sensitivity regions. In particular, two points are near the local maxima of the sensitivity curve, while the other point is placed between them. These design points are seeking to balance the high sensitivity with uniformity and placement of points in regions of greater nonlinearity.

\subsubsection*{Improved trajectory and controller}

\begin{figure}[h]
\centering
  \includegraphics[width=0.4\textwidth]{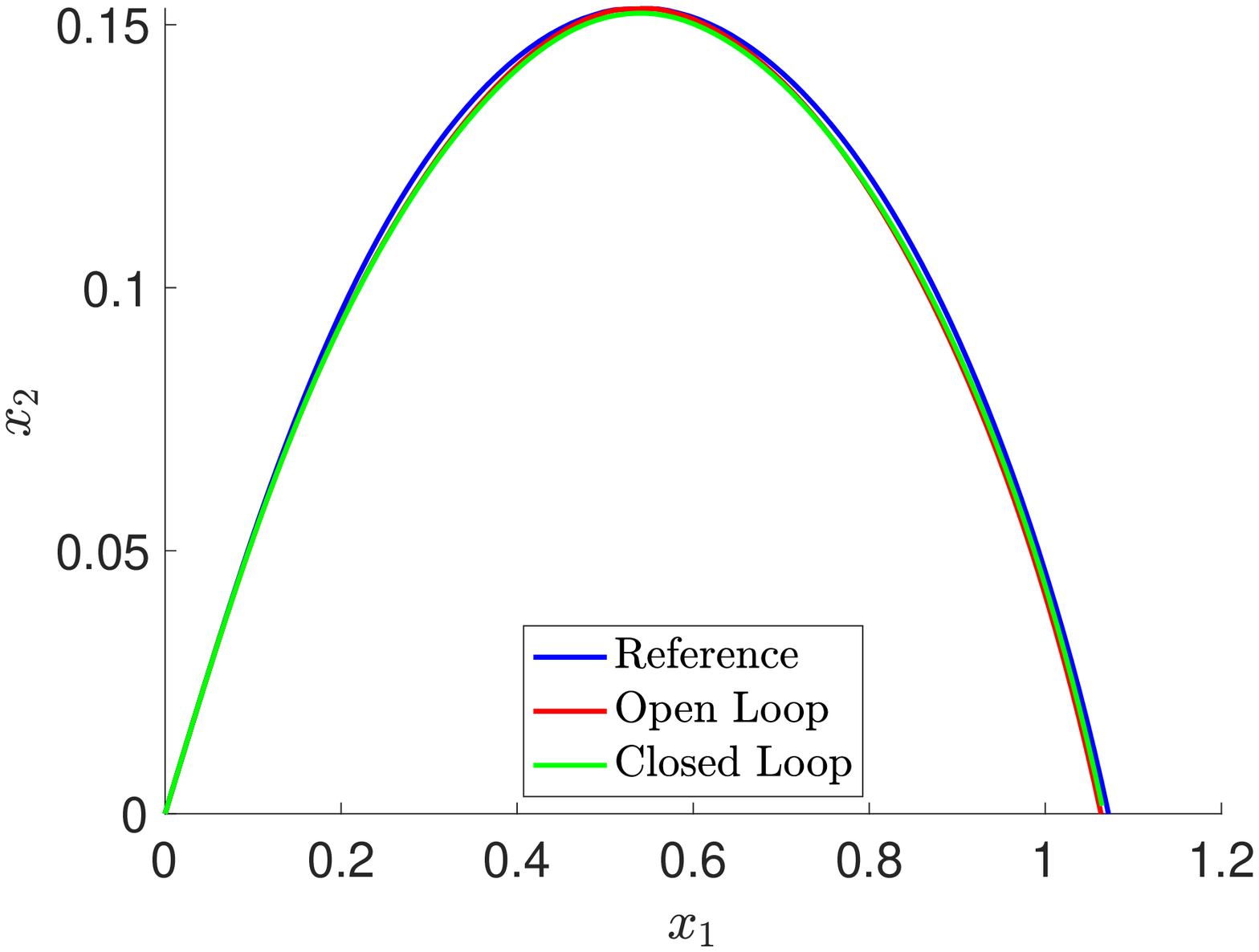}
    \includegraphics[width=0.4\textwidth]{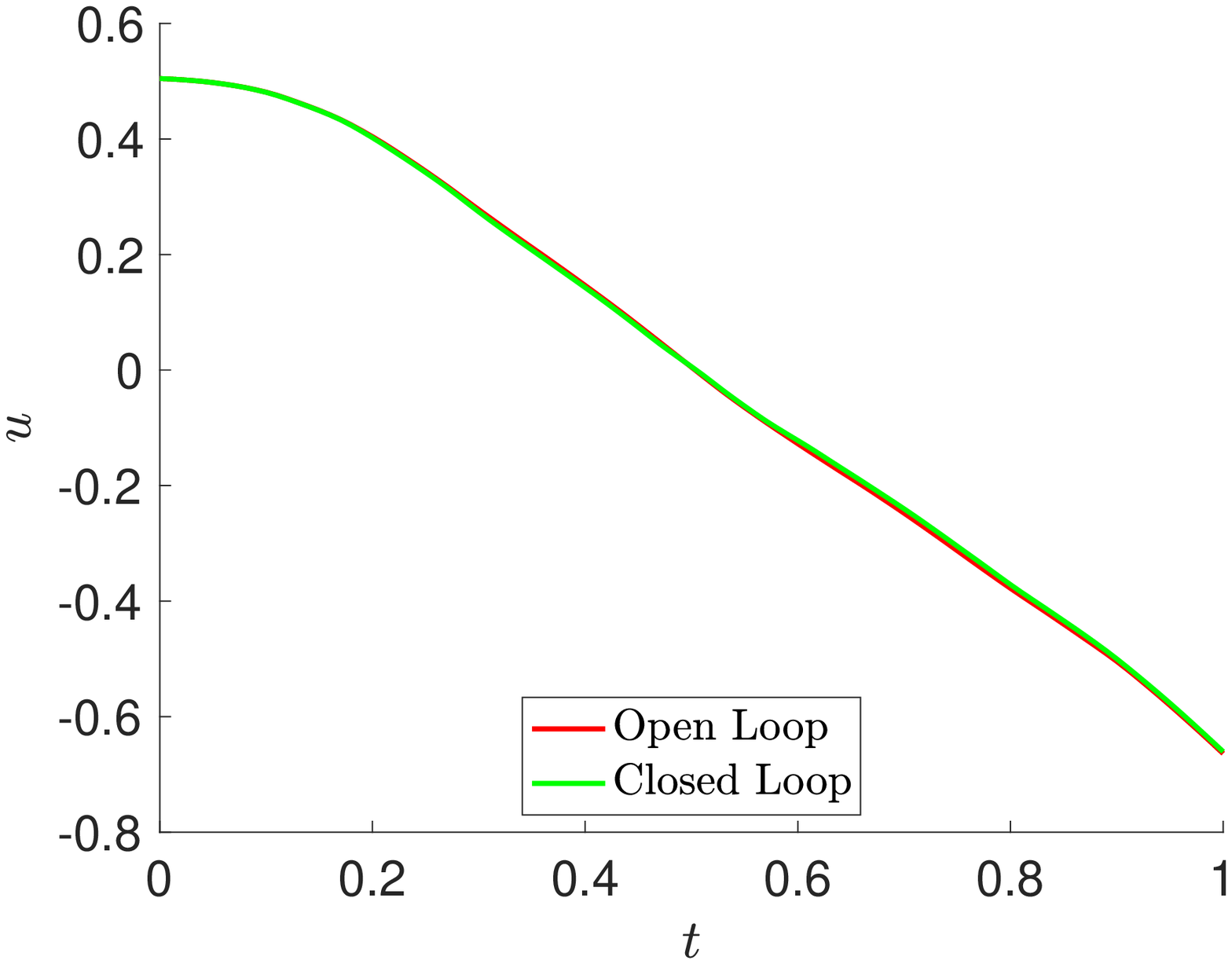}
  \caption{Trajectories in the $x_1 x_2$ phase space (left) and controllers in the time domain (right) with the reference solution $(\tilde{\x},\tilde{u})$ generated using the nominal current model $\tilde{g}$. The reference trajectory $\tilde{\x}$ is the solution of~\eqref{eq:zermelo_OC} with $g=\tilde{g}$ while the open and closed loop trajectories are generated using the controllers $\tilde{u}$ and $\tilde{u}+\Delta u$ with the current model $g^\star$.}
  \label{fig:improved_trajectory}
\end{figure}

We resolve the open loop and LQR feedback control problems using the improved model $g=\tilde{g}$. The updated reference, open loop, and closed loop trajectories are shown in Figure~\ref{fig:improved_trajectory}. We observe that the open loop solution is closer to the reference since the uncertainty in $g$ was reduced. The closed loop solution improves the open solution with a small feedback effort to track the reference trajectory (compare the open and closed loop controllers in the right panel of Figure~\ref{fig:improved_trajectory}). In particular, there is a notable decrease in the feedback control effort from Figure~\ref{fig:trajectory_1} to Figure~\ref{fig:improved_trajectory}. This was the goal of the optimal experimental design.

\subsubsection*{Comparison of experimental designs}
Our previous analysis does not demonstrate how beneficial the proposed experimental design is relative to other options (for instance, an analyst manually choosing where to evaluate $g^\star$ without aid from computation). To address this question we perform a computational test in which we sample different experimental designs. In particular, we generate 100 different random $3$ point designs, and for each design measure the reference tracking error and closed loop control effort. To ensure a fair comparison, for each random design we (1) compute $\tilde{g}$ by augmenting $\overline{g}$ with the new samples, (2) resolve the open loop problem to generate a reference trajectory corresponding to $\tilde{g}$, (3) and determine an LQR feedback with $g=\tilde{g}$. We repeat this experiment for 100 different high fidelity models (generated by perturbing $\overline{g}$ with a random linear combination polynomial basis functions) and average the resulting closed loop control effort and reference tracking error.

\begin{figure}[h]
\centering
  \includegraphics[width=0.4\textwidth]{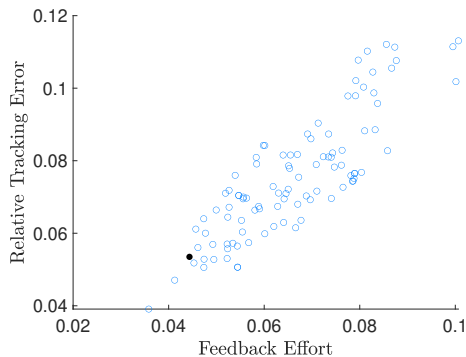}
  \caption{Comparison of experimental designs. Each circle denotes a random design and the dot denotes the optimal experimental design, each averaged over 100 high fidelity models; the horizontal axis indicate the closed loop control effort while the vertical axis is the reference tracking error. An ideal design (which tracks the reference with minimal feedback effort) will be as close to $(0,0)$ as possible.}
  \label{fig:closed_loop_compare}
\end{figure}

Figure~\ref{fig:closed_loop_compare} displays the average feedback control effort and reference tracking error using the 100 random designs as well as the optimal design from Figure~\ref{fig:closed_loop_compare}. 

We observe that the experimental design generated by solving~\eqref{eq:oed_opt_3} is not optimal in the sense of attaining a minimum feedback effort to track the reference trajectory; however, we should not expect it to be. The design generated by solving~\eqref{eq:oed_opt_3} is optimal for the error model specification, information gain assumptions ($r_{i,j}$), and the linear approximation HDSA makes to the parameter to optimal controller mapping. The error and information gain quantities will never be fully known in practice; nonetheless, our framework provides a pragmatic way to incorporate best estimates of them with computational insights in a systemic framework which outperforms the majority of random designs. The design may be easily improved when domain expertise is available, but is not dependent on it.

\subsection{X43 trajectory control} \label{ssec:X43}
Control, guidance, and navigation of hypersonic vehicles is a challenging problem and active research area due to complex nonlinear physics and fast time scales. We demonstrate the utility of our proposed method to control a three degree of freedom point-mass model of a hypersonic glide vehicle using an aerodynamic surrogate model trained on data from the NASA X43 flight tests~\cite{hyperX,x43_2006,x43_reconstruction,parish_model_fidelity_for_RTG}. 

Let $(x,y,z)$ be the position of the vehicle modeled in the north, east, down coordinates, $v$ be the velocity, $\gamma$ be the flight path angle, and $\psi$ be the azimuth angle. The vehicle is controlled by the angle-of-attack $\alpha$, side slip angle $\beta$, and roll angle $\phi$. This gives state and control vectors $\x=(x,y,z,v,\gamma,\psi)^T \in \r^6$ and $\u=(\alpha,\beta,\phi)^T \in \r^3$, respectively. Given an initial state $\x_0 \in \r^6$ and target position $(x_{\text{target}},y_{\text{target}},z_{\text{target}})$, we consider the free final time optimal control problem
\begin{align}\label{eq:hypersonic_OC}
& \min_{\x,\u,T}  (x-x_{\text{target}})^2 + (y-y_{\text{target}})^2 + (z-z_{\text{target}})^2\\
& s.t. \nonumber \\
& \begin{dcases} \label{eqn:traj_opt}
 \dot{x} = v\cos(\gamma) \cos(\psi)  & t \in (0,T) \\
 \dot{y} = v\cos(\gamma) \sin(\psi)  & t \in (0,T) \\
 \dot{z} = - v \sin(\gamma)  & t \in (0,T) \\
 \dot{v} = - \frac{D}{m}- g \sin(\gamma)& t \in (0,T)  \\
 \dot{\gamma} = \frac{L \cos(\phi) - S \sin(\phi)}{m v} - \frac{g \cos(\gamma)}{v}  \qquad  & t \in (0,T)  \\
 \dot{\psi} = \frac{ L \sin(\phi) + S \cos(\phi)}{mv\cos(\gamma)} & t \in (0,T)  \\
 \x(0) = \x_0
 \end{dcases}\\
 & \begin{dcases} \label{eqn:x43_inequal}
k_Q \sqrt{\rho} v^3 \le Q_{\text{max}} \qquad & t \in (0,T) \\
\frac{1}{2} \rho v^2 \le q_{\text{max}} \qquad & t \in (0,T) \\
\frac{\sqrt{L^2 + D^2}}{m} \le n_{\text{max}} \qquad & t \in (0,T) \\
(x-c_x)^2 + (y-c_y)^2 \le r_c^2 \qquad & t \in (0,T) \\
\alpha_l \le \alpha \le \alpha_u \qquad & t \in (0,T) \\
\beta_l \le \beta \le \beta_u \qquad & t \in (0,T) \\
\phi_l \le \phi \le \phi_u \qquad & t \in (0,T) \\
 \end{dcases}
\end{align}
where~\eqref{eqn:x43_inequal} are path constraints corresponding to the maximum heating rate $Q_{\text{max}}$, dynamic pressure $q_{\text{max}}$, load factor $n_{\text{max}}$, and a cylindrical no flight zone centered at $(x_c,y_c)$ with radius $r_c$, and controller bounds $\alpha_l,\alpha_u, \beta_l,\beta_u,\phi_l,\phi_u$.

The ODE system~\eqref{eqn:traj_opt} and path constraints~\eqref{eqn:x43_inequal} depend on the aerodynamics model for drag, lift, and side forces
\begin{align*}
&D  = (F_N \sin(\alpha) + F_A \cos(\alpha)) \cos(\beta) - F_Y \sin(\beta) \\
&L  = F_N \cos(\alpha) - F_A \sin(\alpha) \\
&S  = (F_N \sin(\alpha) + F_A \cos(\alpha)) \sin(\beta) + F_Y \cos(\beta) \\
\end{align*}
where
\begin{align*}
&F_N  = \frac{1}{2} \rho v^2 A_{ref} C_N \\
&F_A  = \frac{1}{2} \rho v^2 A_{ref} C_A \\
&F_Y  = \frac{1}{2} \rho v^2 A_{ref} C_Y
\end{align*}
are the aerodynamic forces in the normal, axial, and yaw directions and $C_N(M,\alpha)$, $C_A(M,\alpha)$, and $C_Y(M,\alpha,\beta)$ are the non-dimensional aerodynamic coefficients which depend on the Mach number (a function of velocity and altitude) $M$, angle-of-attack $\alpha$, and in the case of $C_Y$, the side slip angle $\beta$.  An exponential model $\rho=\rho_0 \exp(z/H)$ is used to emulate the atmospheric density. The constants $m$, $A_{ref}$, and $k_Q$ are the mass, reference area, and heat flux coefficient of the vehicle, respectively; $g$ is the acceleration due to gravity. All results will be reported in normalized units.

 We consider uncertainty in the aerodynamic coefficients
$$\g(M,\alpha,\beta)=(C_N(M,\alpha),C_A(M,\alpha),C_Y(M,\alpha,\beta))^T \in \r^3$$
which are determined by fitting a surrogate model to aerodynamic data. Notice that $\g$ depends on $M$ which is a function of state variables rather than a state variable itself, and that the third component of $\g$ depends on an additional variable, $\beta$, which the first two components do not. 
   
  \subsubsection*{Trajectory and controller}
A nominal estimate for the aerodynamic coefficients, $\gbar$, is given by a polynomial model constructed from the NASA X43 flight data. To illustrate the concepts of this article, we define the true aerodynamics $\gstar$, which are not available in practice, as a perturbation of $\gbar$ which we defined so that we can analyze the controller performance (perform experiments with $\gstar$). We solve the trajectory optimization problem~\eqref{eq:hypersonic_OC} with $\overline{\g}$ to generate a reference trajectory and controller $(\overline{\x},\overline{\u})$ and adopt an LQR feedback controller as in the Zermelo problem. 

\begin{figure}[h]
\centering
          \includegraphics[width=0.5\textwidth]{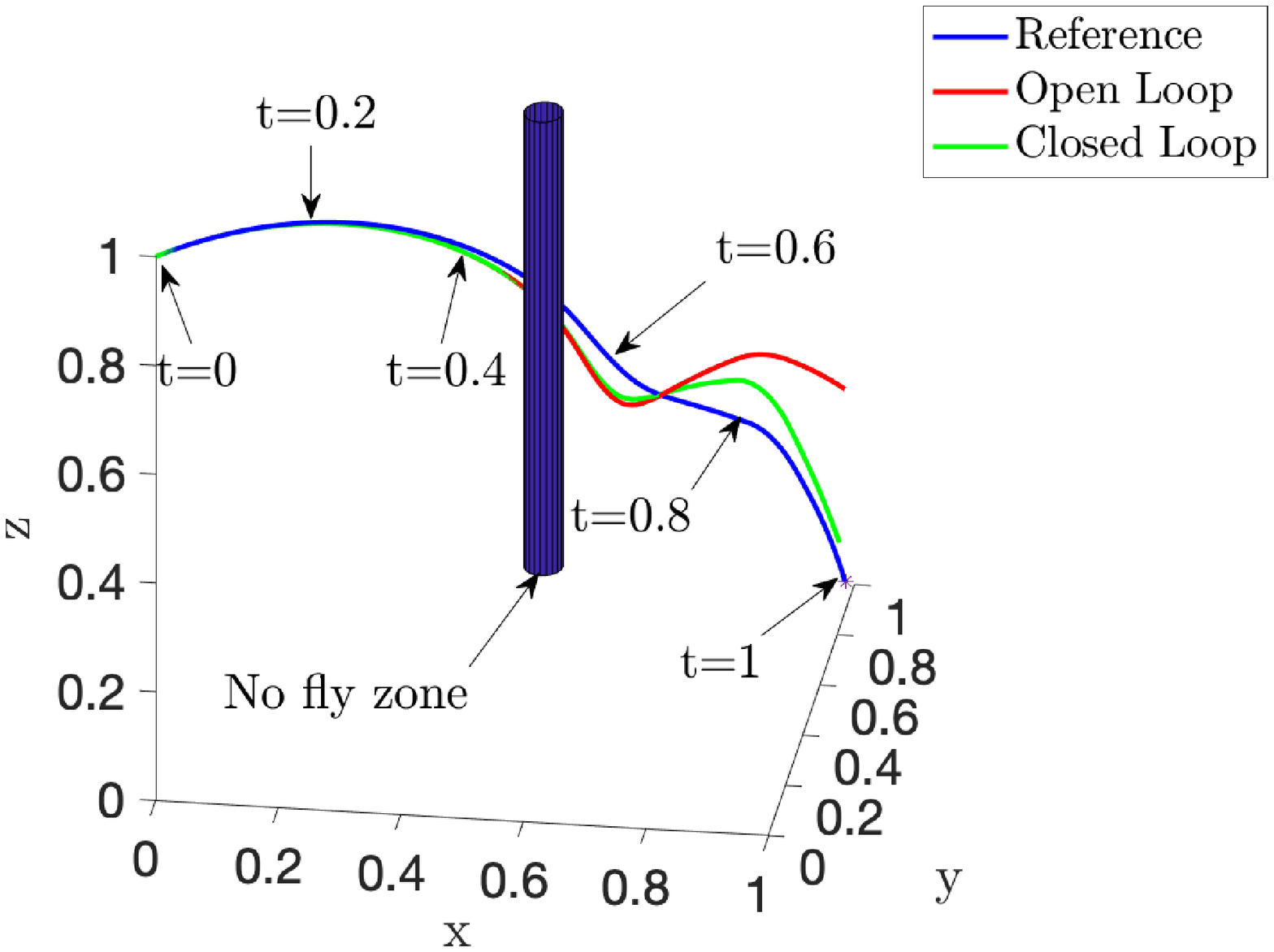}\\
  \includegraphics[width=0.32\textwidth]{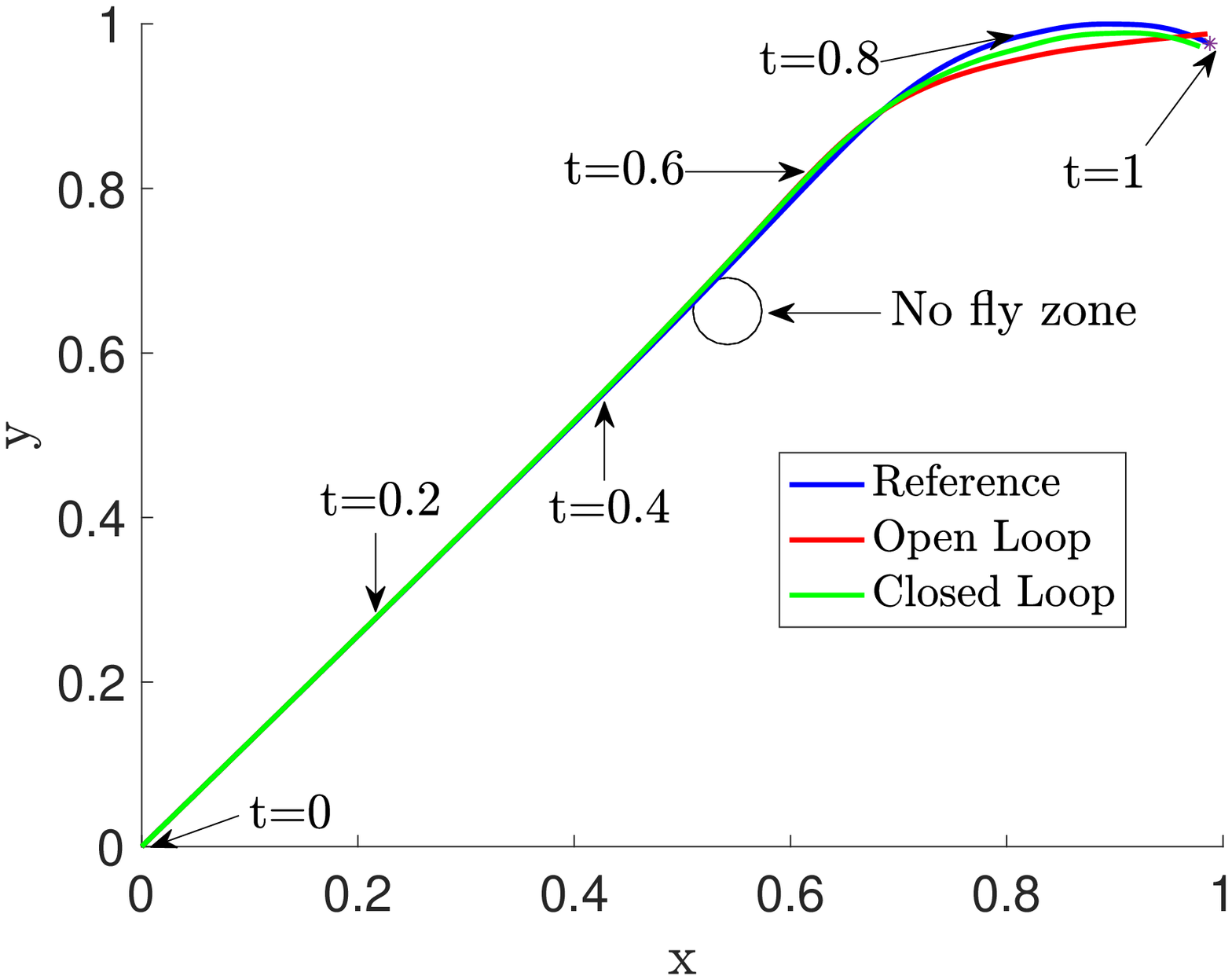}
    \includegraphics[width=0.32\textwidth]{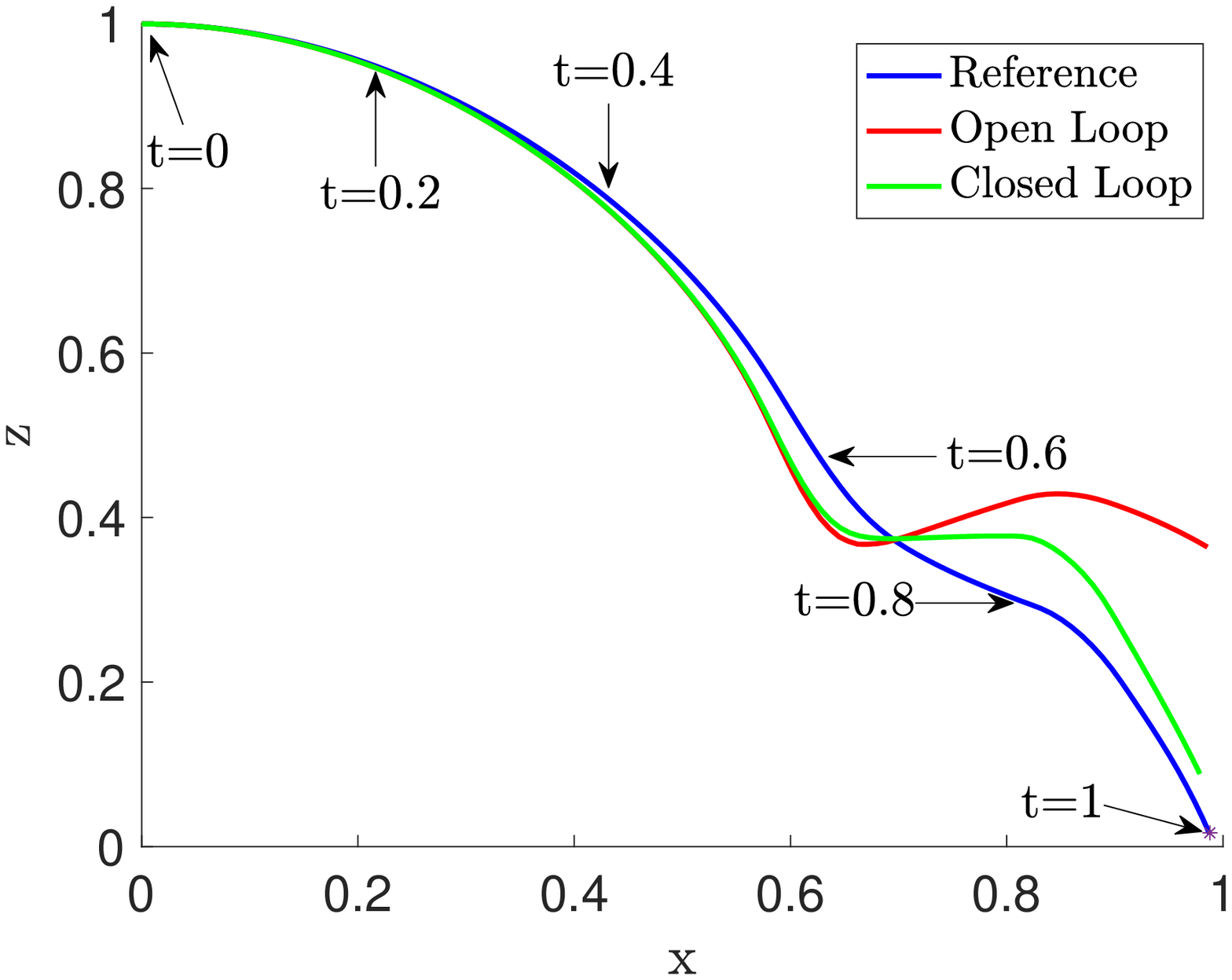}
      \includegraphics[width=0.32\textwidth]{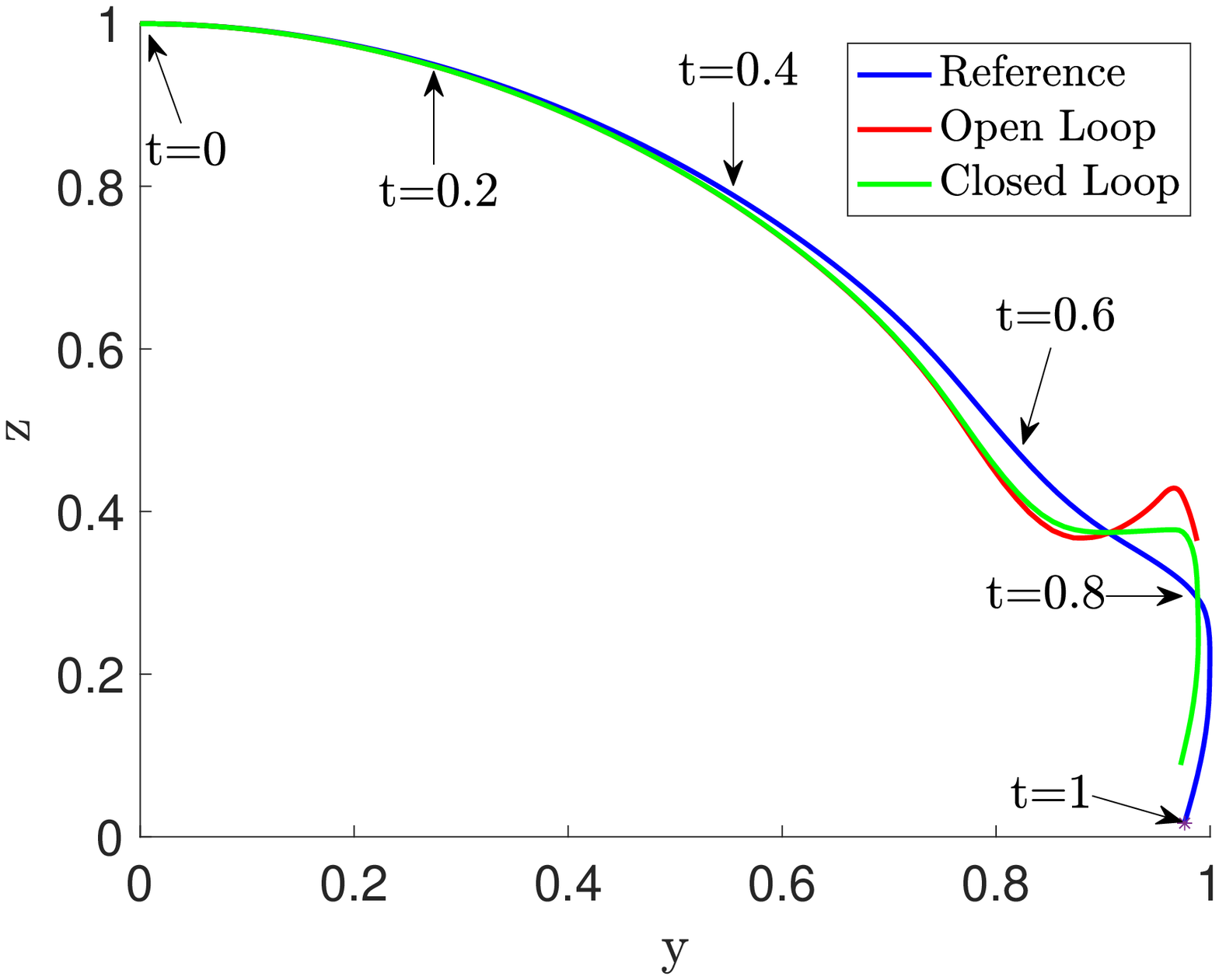}
  \caption{Trajectories with a reference solution generated using the nominal aerodynamic model $\gbar$. The reference trajectory $\xbar$ is the solution of~\eqref{eq:hypersonic_OC} with $\g=\gbar$ while the open and closed loop trajectories are generated using the aerodynamic model $\gstar$. The top panel shows the trajectory in three dimensional space while the bottom row shows two dimensional views of the trajectory for easier visualization.}
  \label{fig:x43_trajectory_1}
\end{figure}

Figure~\ref{fig:x43_trajectory_1} displays the reference trajectory $\overline{\x}$, the open loop trajectory $\xo$, and the closed loop trajectory $\xc$, in the spatial coordinates $(x,y,z)$. In particular, the top panel shows the trajectories in 3D, while the bottom row shows each 2D cross section. The no fly zone is shown by the cylinder and time stamps mark the location of the vehicle at times $t=0,0.2,0.4,0.6,0.8,1.0$. We observe that the open loop trajectory $\xo$ diverges far from the reference trajectory $\xbar$ due to errors in the aerodynamic coefficients. The addition of the LQR feedback yields a closed loop trajectory $\xc$ which is improved, but is still insufficient to track the reference.
  
\begin{figure}[h]
\centering
  \includegraphics[width=0.32\textwidth]{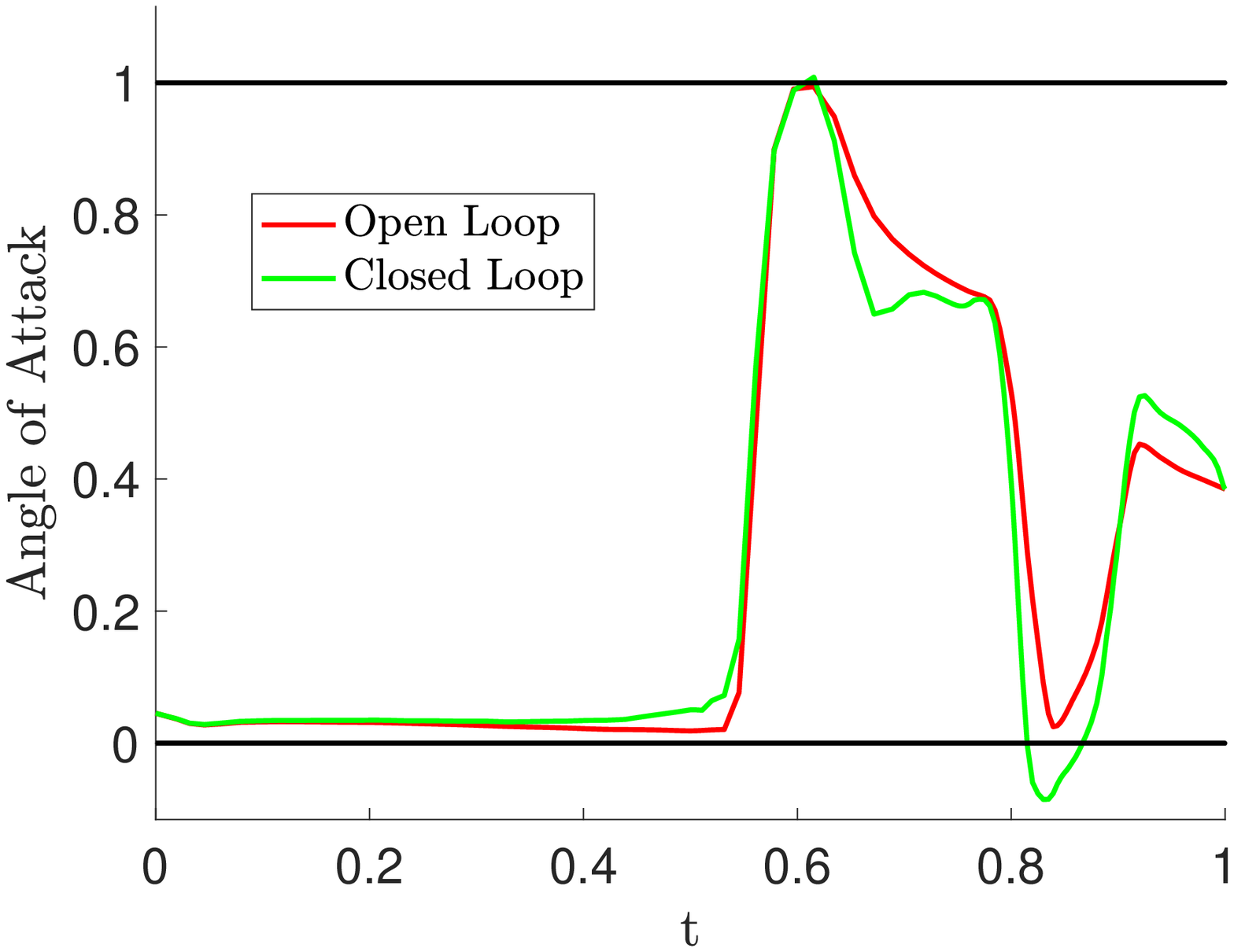}
    \includegraphics[width=0.32\textwidth]{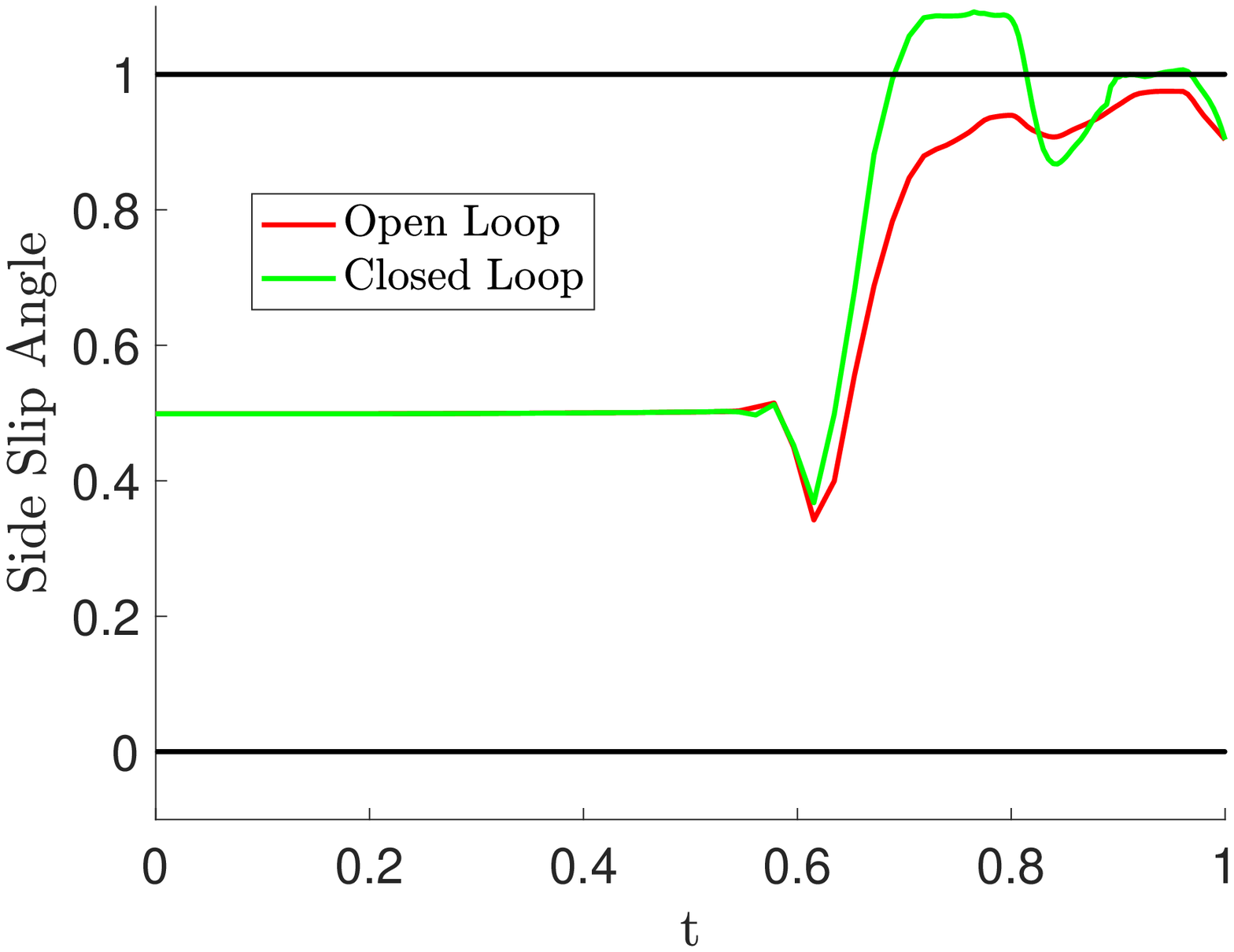}
      \includegraphics[width=0.32\textwidth]{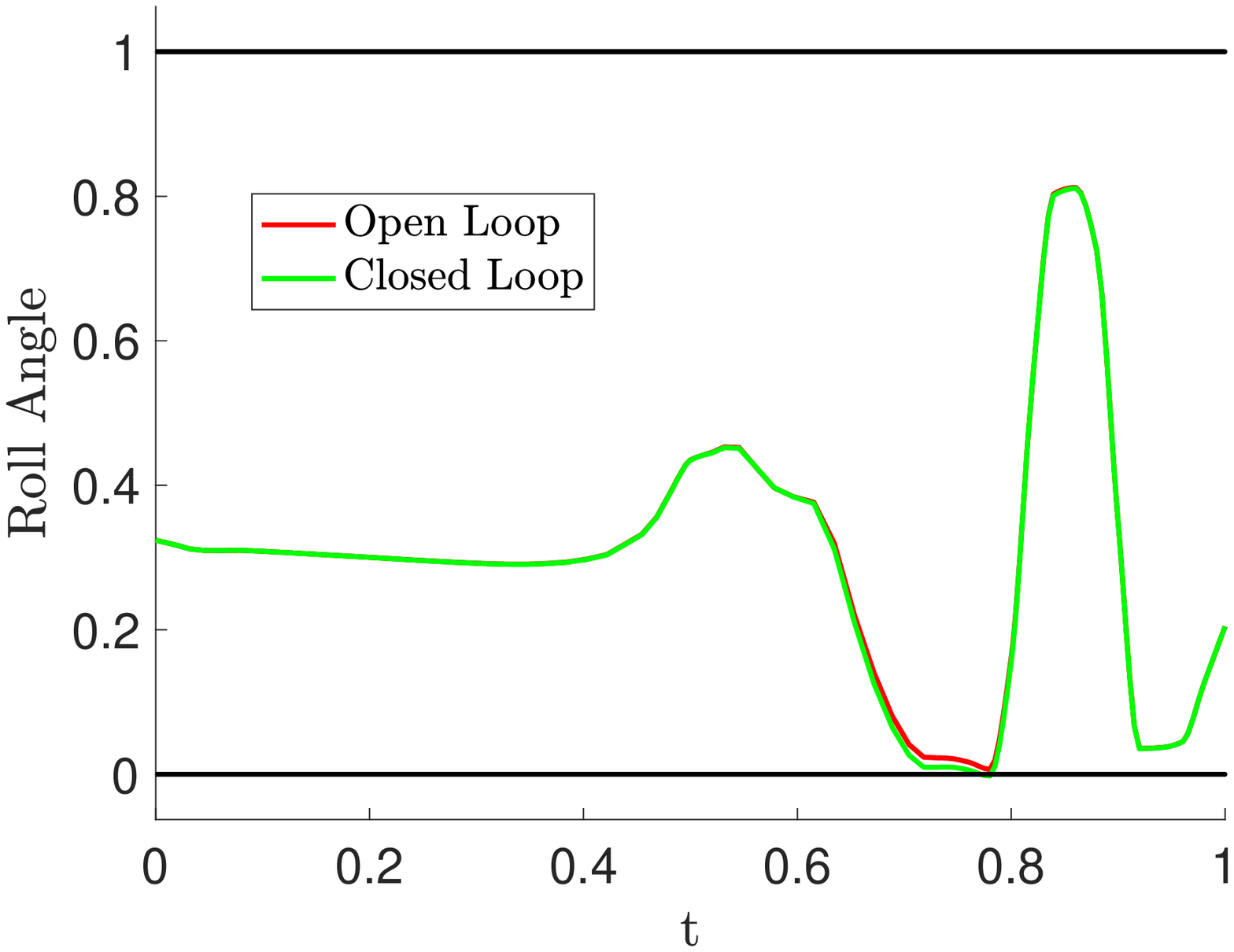}
  \caption{Controllers used to generate the trajectories in Figure~\ref{fig:x43_trajectory_1}. The open loop controller $\ubar$ is generated by solving~\eqref{eq:hypersonic_OC} with $\g=\gbar$ while the closed loop controller $\ubar + \fb$ uses LQR feedback around the nominal trajectory. The horizontal lines denote the lower and upper bounds on the controllers.}
  \label{fig:x43_control_1}
\end{figure}

Figure~\ref{fig:x43_control_1} shows the controllers corresponding to the trajectories in Figure~\ref{fig:x43_trajectory_1}. The open loop controller is $\ubar$, computed by solving~\eqref{eqn:traj_opt}, and the closed loop controller is $\ubar + \fb$, generated using the LQR feedback controller. The horizontal lines correspond to the upper and lower bounds on the controllers~\eqref{eqn:x43_inequal}. Considerable feedback effort is needed after $t=0.6$ and it violates the bounds on the controller. This large feedback effort, which improves the trajectory but still fails to track the reference, is what we seek to reduce via better characterization of the aerodynamics in the most important flight configurations.

\subsubsection*{Hyper-differential sensitivities and experimental design}
\begin{figure}[h]
\centering
   \includegraphics[width=0.5\textwidth]{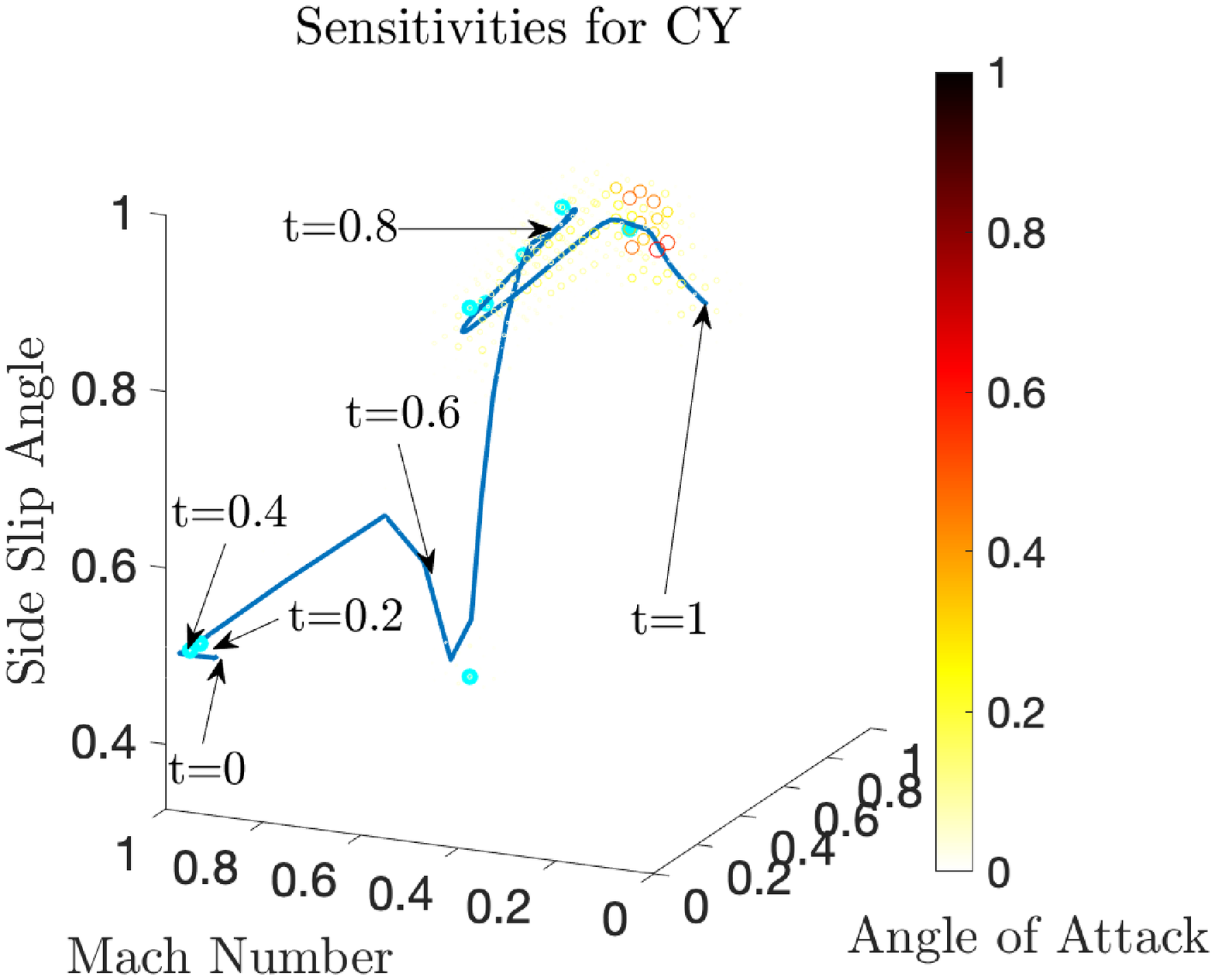}\\
   \includegraphics[width=0.32\textwidth]{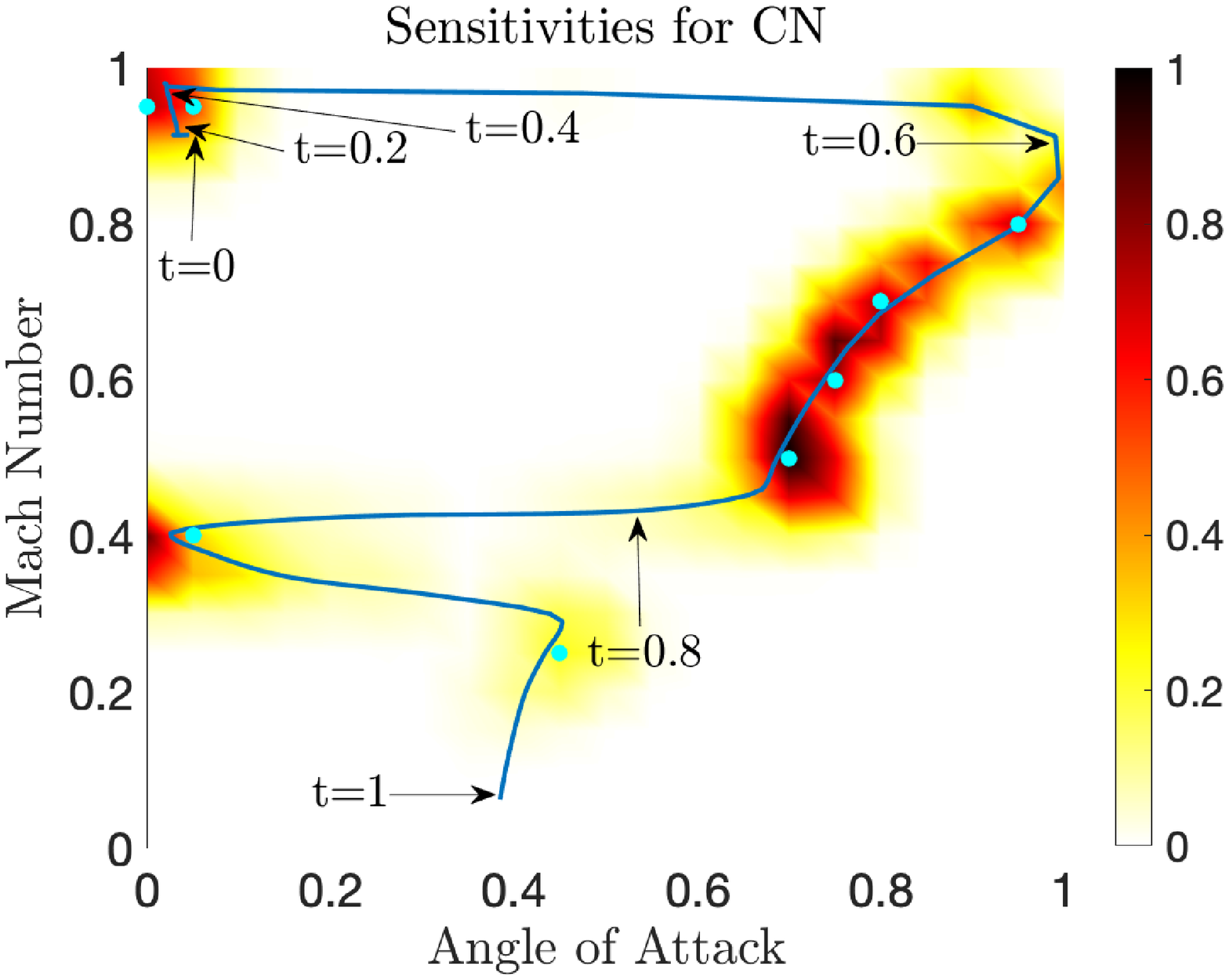}
  \includegraphics[width=0.32\textwidth]{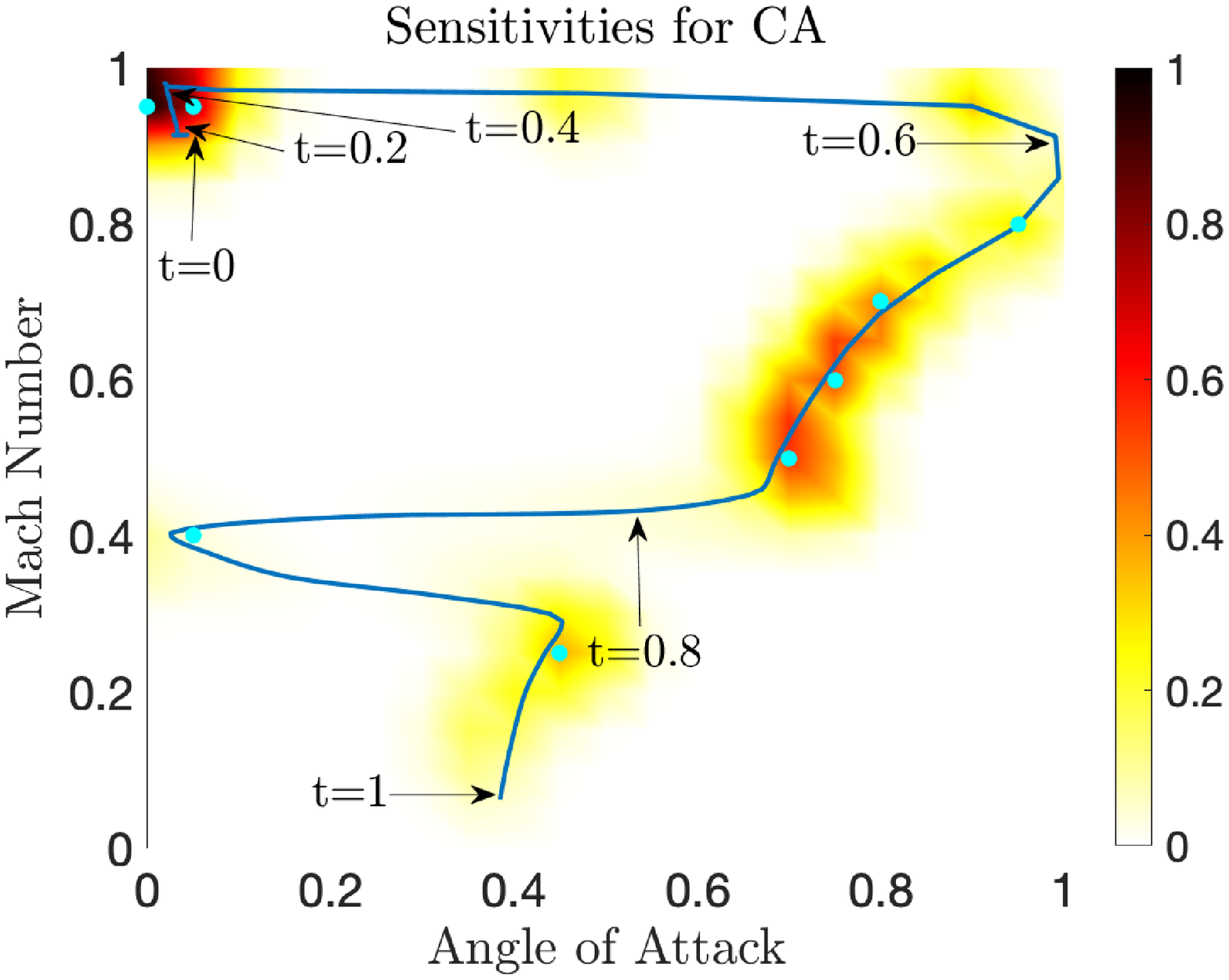}
        \includegraphics[width=0.32\textwidth]{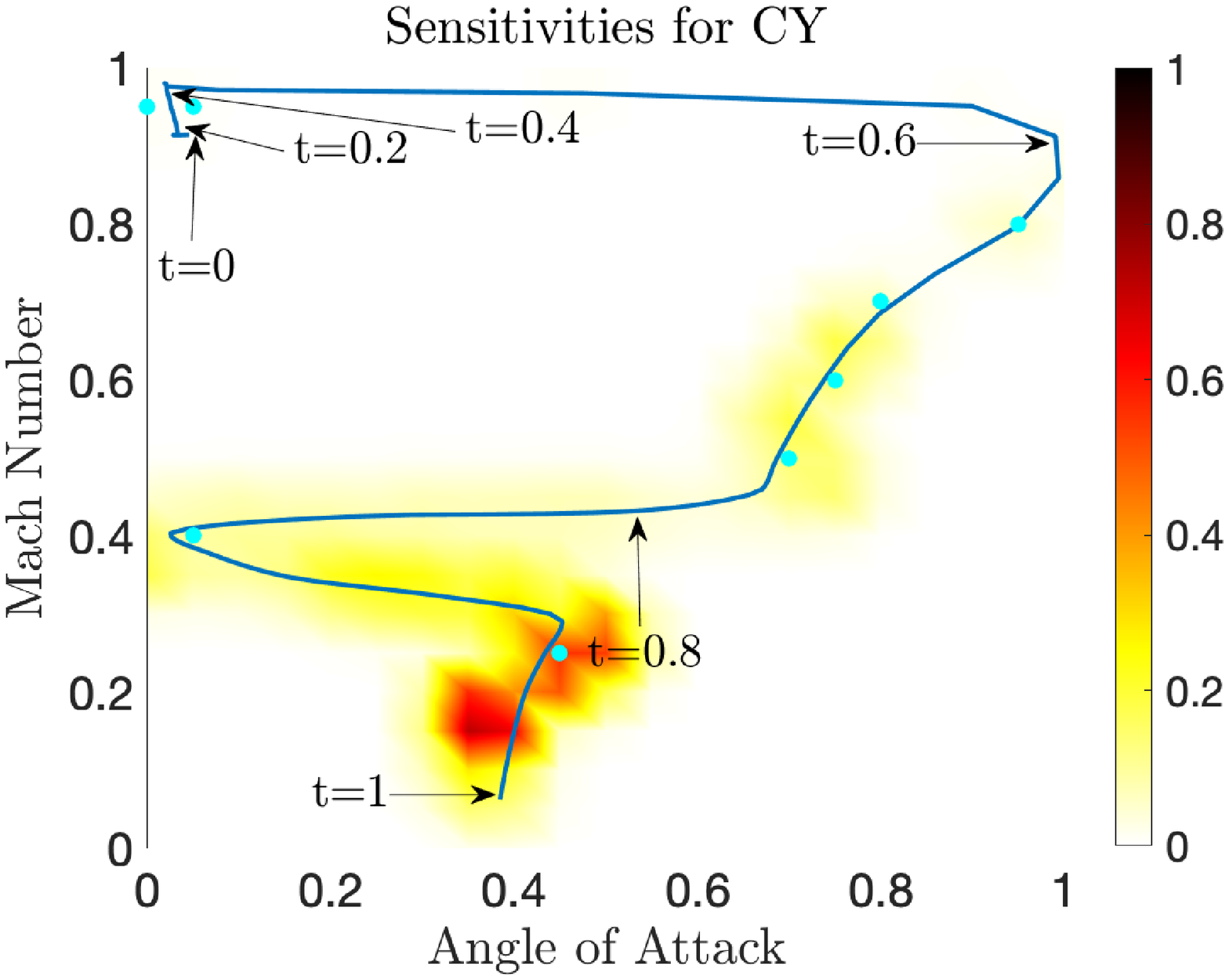}
  \caption{Hyper-differential sensitivities for the X43 reference tracking problem. The top panel is $C_Y$ in its three dimensional phase space for $(M,\alpha,\beta)$, while the bottom row is $C_N$ and $C_A$ in their two dimensional space space for $(M,\alpha)$, along with $C_Y$ projected onto the $(M,\alpha)$ phase space. The blue curve corresponds to the reference trajectory $\xbar$ while the eight cyan dots indicate the location of experimental design points generated using the sensitivities.}
  \label{fig:x43_sensitivities}
\end{figure}

We compute hyper-differential sensitivities to determine which flight configurations $(M,\alpha,\beta)$ result in aerodynamic uncertainty placing the greatest demand on the feedback controller. Figure~\ref{fig:x43_sensitivities} displays the hyper-differential sensitivities, plotted in the $(M,\alpha,\beta)$  phase space, with the reference trajectory $\xbar$ overlaid. Note that the aerodynamic coefficient $C_N$ and $C_A$ only depend on $M$ and $\alpha$ so we plot them in a two dimensional phase space. Since $C_Y$ also depends on $\beta$, we provide two plots to show the sensitivities. The top panel of Figure~\ref{fig:x43_sensitivities} is the three dimensional view while the bottom right panel projects $C_Y$ onto the $(M,\alpha)$ phase space by averaging in the $\beta$ coordinate and renormalizing so that the maximum sensitivity in the $(M,\alpha)$ phase space equals the maximum in the $(M,\alpha,\beta)$ phase space. This additional normalization is necessary because the sensitivities are along a two dimensional curve in a three dimensional phase space and hence typical averaging is not appropriate.
 
We observe that the greatest sensitivity corresponds to axial force ($C_A$) during the initial glide phase $t \in (0,.5)$ and the normal force ($C_N$) during the diving phase $t \in (0.65,0.75)$. This is because the axial force during the glide phase dictates the range of the vehicle over this longer time horizon while the normal force during the dive phase determines how quickly the vehicle descends. Notice that the high axial force sensitivity during the glide phase, $t \in (0,.5)$, results in needing feedback during the dive phase. This is because the vehicle is less controllable at higher altitude due to the smaller atmospheric density, so the range error induced by poor estimation of aerodynamics in this portion of the trajectory must be compensated for later in the flight. The high normal force sensitivity during the diving phase, $t \in (0.65,0.75)$, also demands feedback during the same time interval. The cumulative effect of these uncertainties explains why considerable feedback effort is needed for $t>0.6$.

An optimal experimental design is computed to determine $\kappa_B = 8$ points in the $(M,\alpha,\beta)$ phase space were we evaluate the true aerodynamics $\gstar$. These points are indicated by the cyan dots in Figure~\ref{fig:x43_sensitivities} which are concentrated in the high sensitivity region and spread in such a way as to optimize the information gain by interpolating between them. In this example, we define the error reduction model $r_{i,j}$ using the derivatives of $\gbar$ so that it reflects smoothness in the aerodynamics. 

\subsubsection*{Improved trajectory and controller}

\begin{figure}[h]
\centering
  \includegraphics[width=0.5\textwidth]{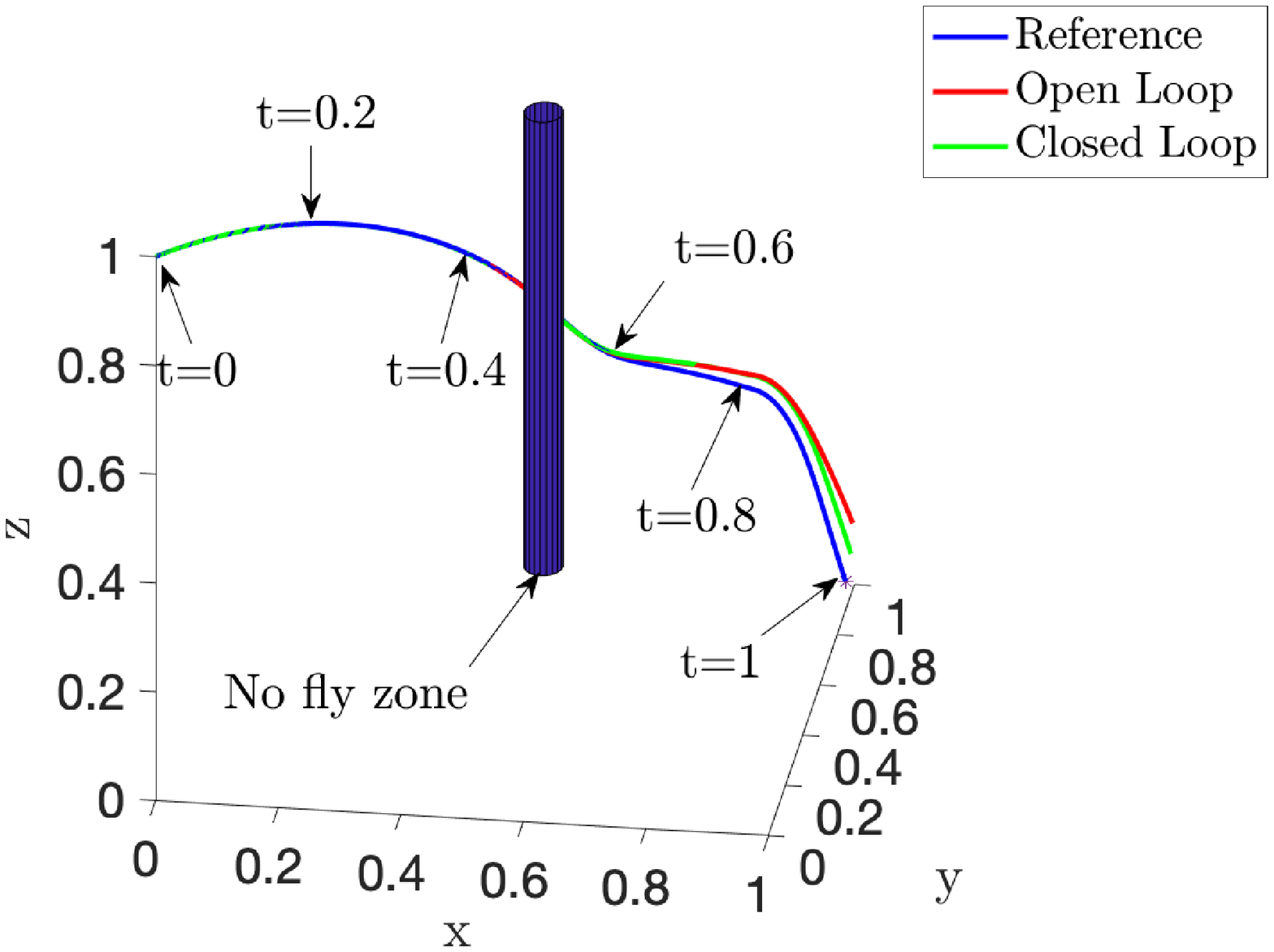}\\
  \includegraphics[width=0.32\textwidth]{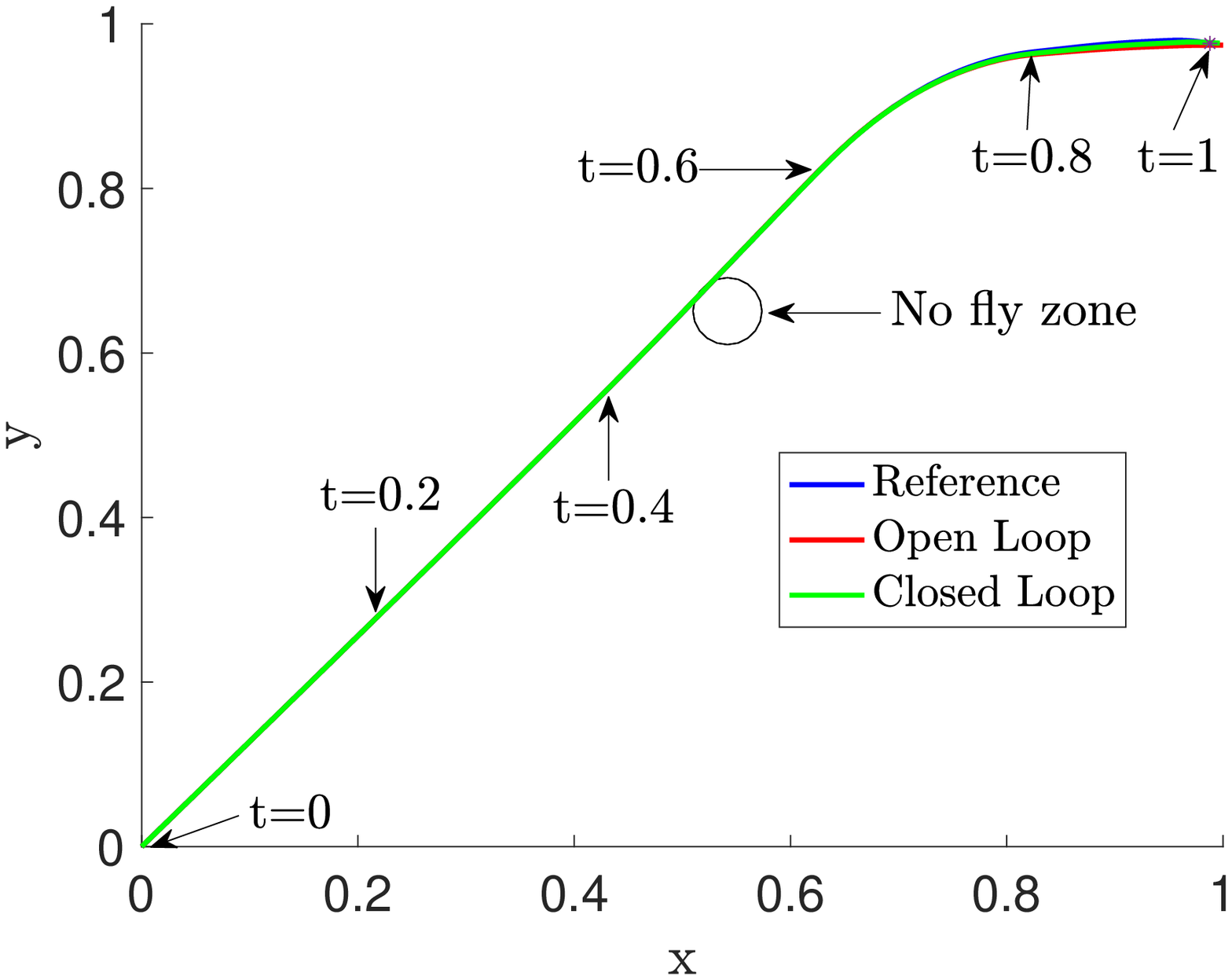}
    \includegraphics[width=0.32\textwidth]{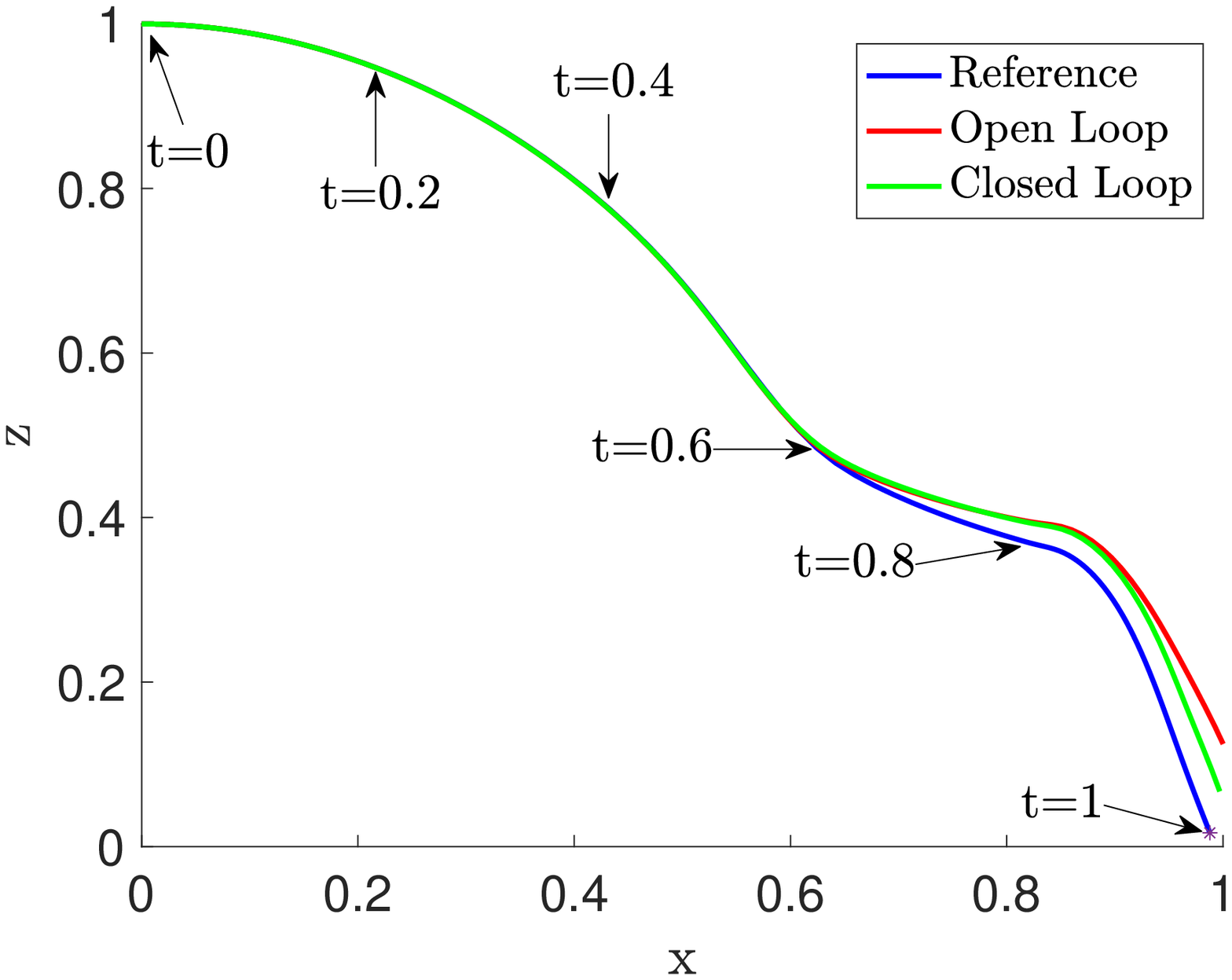}
      \includegraphics[width=0.32\textwidth]{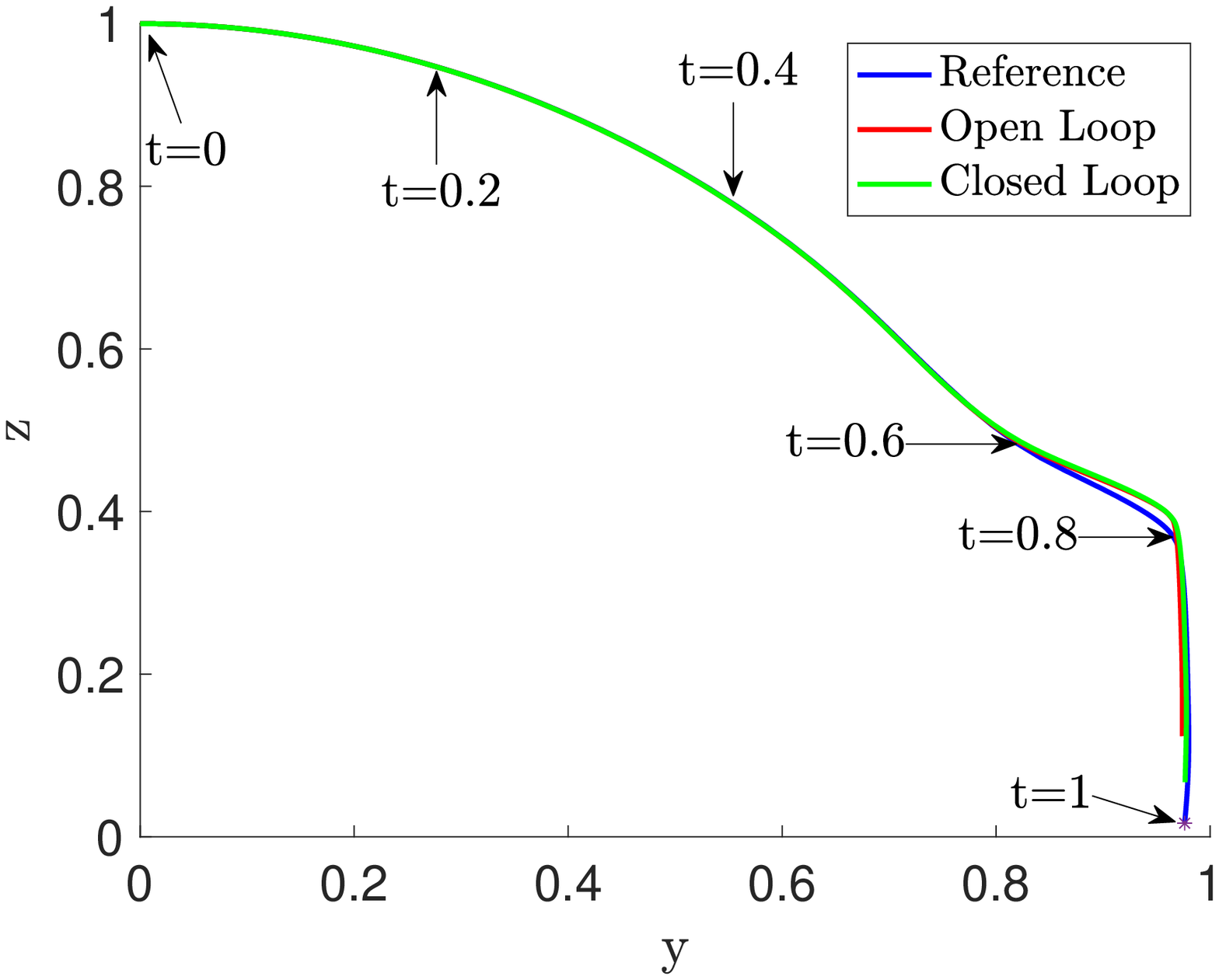}
    \caption{Trajectories with a reference solution generated using the nominal aerodynamic model $\tilde{\g}$. The reference trajectory $\tilde{\x}$ is the solution of~\eqref{eq:hypersonic_OC} with $\g=\tilde{\g}$ while the open and closed loop trajectories are generated using the aerodynamic model $\gstar$. The top panel shows the trajectory in three dimensional space while the bottom row shows two dimensional views of the trajectory for easier visualization.}
  \label{fig:x43_trajectory_2}
\end{figure}

\begin{figure}[h]
\centering
  \includegraphics[width=0.32\textwidth]{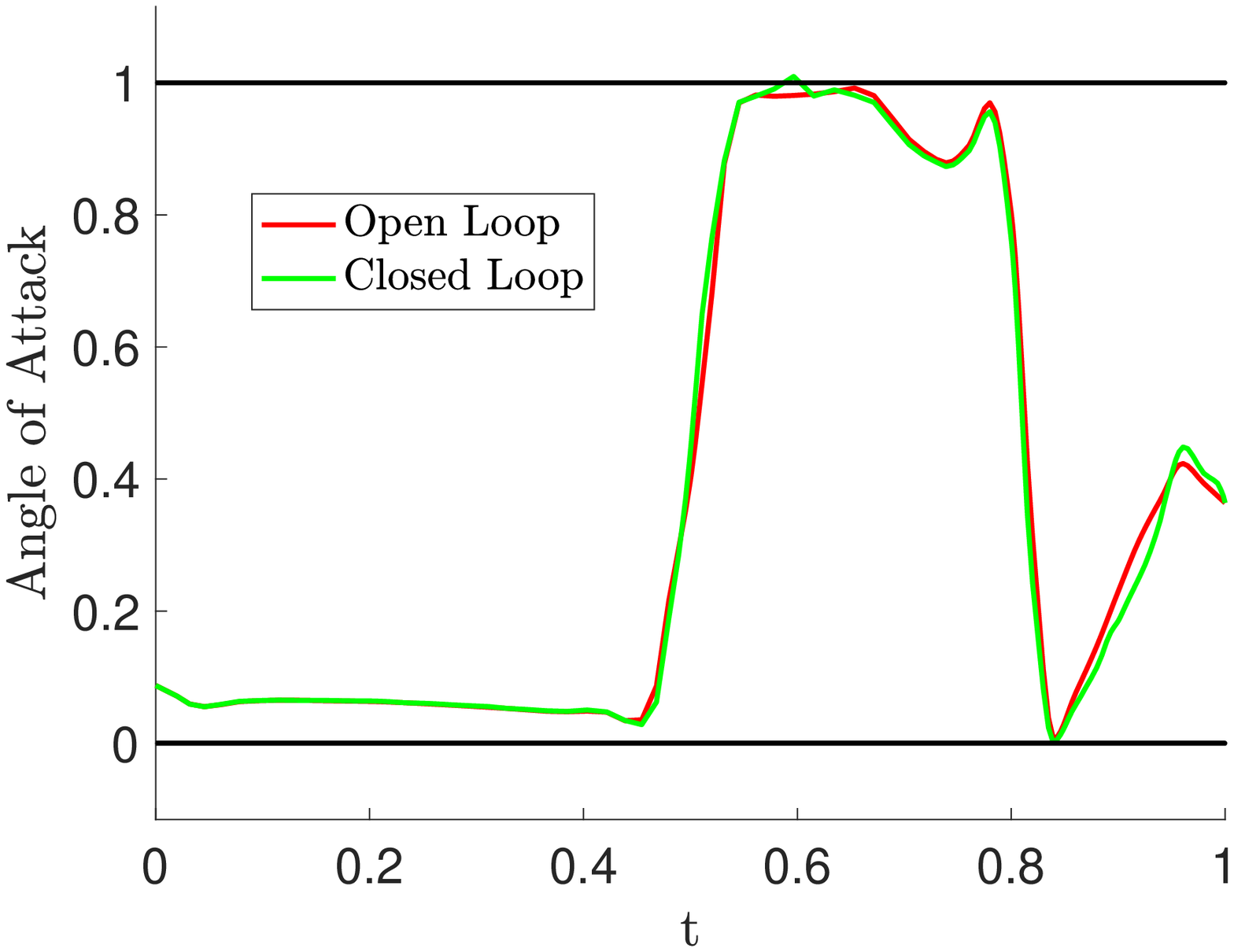}
    \includegraphics[width=0.32\textwidth]{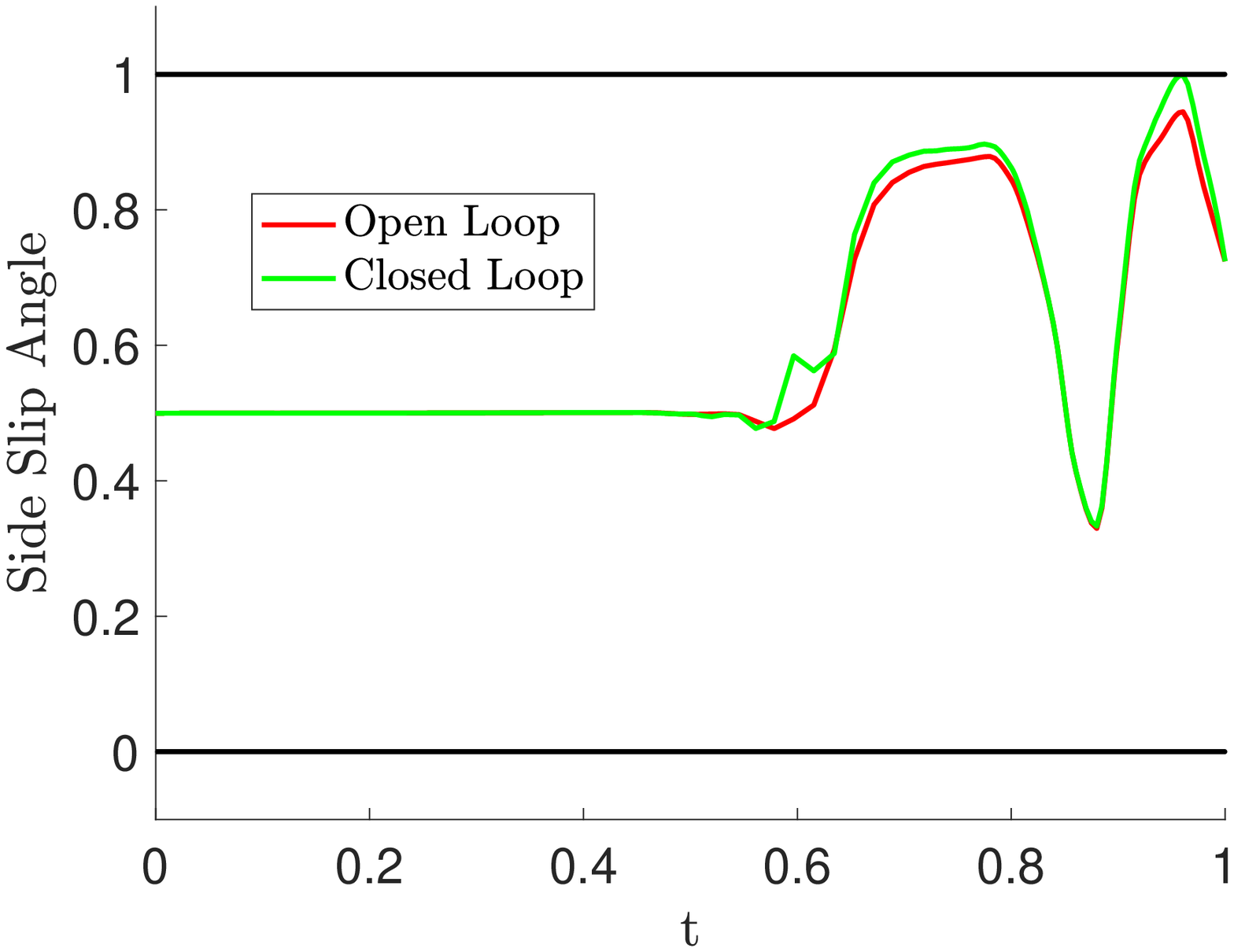}
      \includegraphics[width=0.32\textwidth]{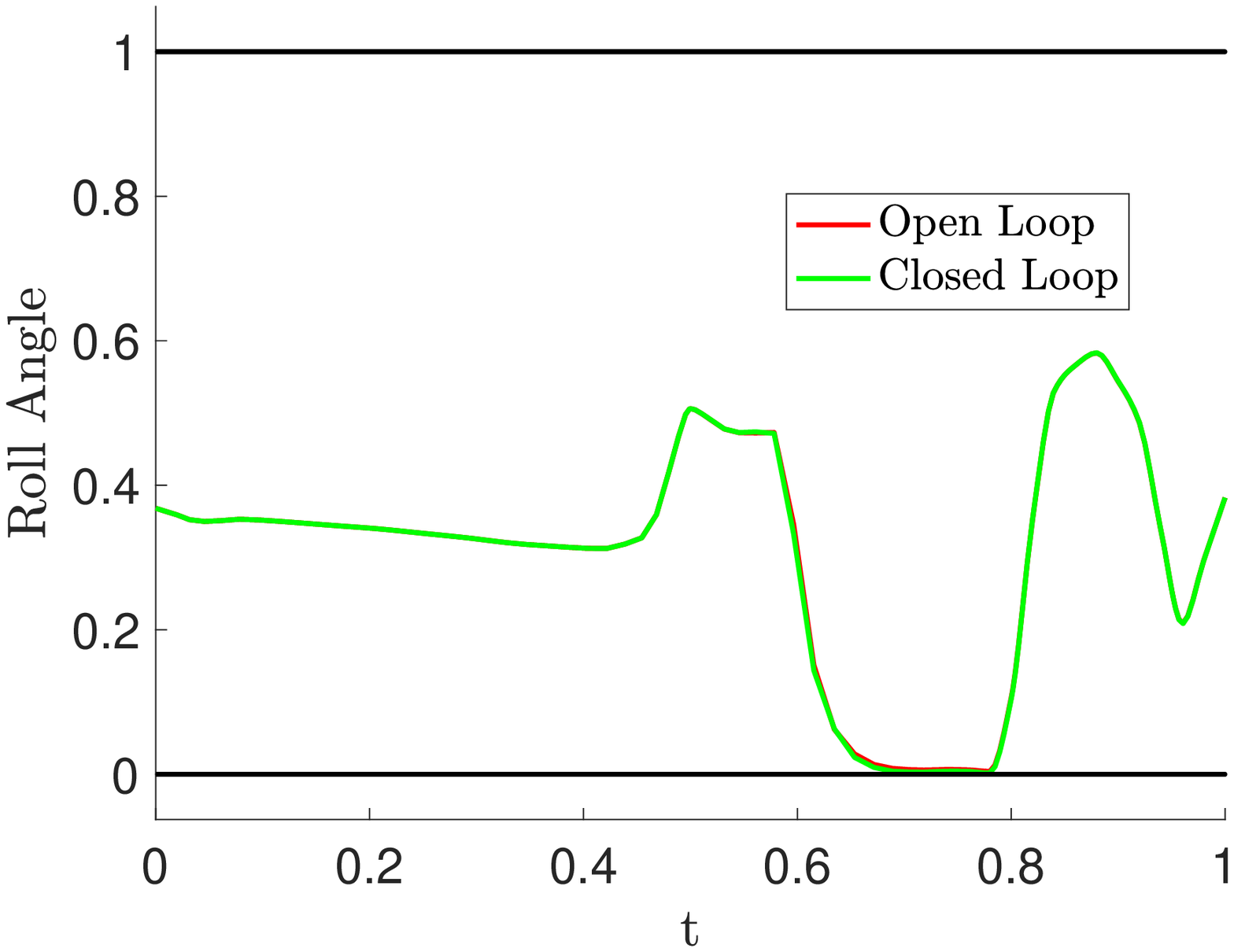}
  \caption{Controllers used to generate the trajectories in Figure~\ref{fig:x43_trajectory_2}. The open loop controller $\tilde{\u}$ is generated by solving~\eqref{eq:hypersonic_OC} with $\g=\tilde{\g}$ while the closed loop controller $\tilde{\u}+\fb$ uses LQR feedback around the nominal trajectory. The horizontal lines denote the lower and upper bounds on the controllers.}
  \label{fig:x43_control_2}
\end{figure}

Using the optimal experimental design, we evaluate the high fidelity model $\gstar$ at $\kappa_B = 8$ points and use this data to improve the aerodynamic model. As in the Zermelo problem, $\tilde{\g}$ is defined by adding radial basis functions to $\gbar$ so that $\tilde{\g}$ interpolates $\gstar$ at the design points. The trajectory planning problem~\eqref{eq:hypersonic_OC} is resolved with the improved aerodynamic model $\tilde{\g}$. Figure~\ref{fig:x43_trajectory_2} shows the reference trajectory $\tilde{\x}$, open loop trajectory, and closed loop trajectory (generated using an LQR feedback) with this improved aerodynamic model. We observe that the closed loop trajectory does not perfectly track the reference $\tilde{\x}$; however, it is a significant improvement relative to the previous trajectory in Figure~\ref{fig:x43_trajectory_1}. More notably, Figure~\ref{fig:x43_control_2} shows the open and closed loop controllers corresponding to these trajectories. We observe that the feedback effort needed in the closed loop setting is considerably less than what was previously required in Figure~\ref{fig:x43_control_1}. This illustrates the utility of the proposed approach in that we were able to avoid violating the control bounds and significantly reduce the demand on the feedback controller.

\section{Conclusion}\label{sec:conclusion}
Optimal control in the face of uncertainty and nonlinearity is challenging in many respects. A variety of approaches exists for both generating nominal trajectories which are robust to uncertainty and designing feedback controllers which aid to mitigate uncertainties. The proposed hyper-differential sensitivity analysis (HDSA) and optimal experimental design approach introduced in this article complements these approaches by seeking to prioritize uncertainties and direct data acquisition which may reduce them. For challenging problems such as autonomous hypersonic flight, wind tunnel testing, computational fluid dynamics, trajectory planning, and feedback control must be organically coupled to overcome the vast uncertainty, strong nonlinearity, and fast time scales faced in practice. As a tool which bridges this areas, hyper-differential sensitivity analysis is poised to contribute invaluable insight and guidance for vehicle development, flight planning, and onboard control. The interplay between uncertainty calibration, trajectory planning, and feedback control pervades many applications and hence the proposed approach, which was demonstrated for control of a hypersonic vehicle, has broad applicability in the engineering sciences.

The contributions of this article focused on how HDSA directs data acquisition in the service of improving trajectory planning. Because of its interpretation as the feedback control effort required to overcome uncertainty, there is opportunity to leverage HDSA for other aspects of both trajectory planning and feedback controller design. By analyzing the sensitivity of the reference tracking optimization problem, rather than simply the model dynamics, HDSA provides unique insights which support decision-making.

Our sensitivity driven optimal experimental design presents a new perspective on data acquisition for optimal control problems. Traditionally, experimental design has focused on statistical estimation and seeks to achieve desirable statistical properties. Our proposed approach is deterministic and grounded in optimal control, yet mimics many important characteristics of its statistical counterpart. 

\section*{Acknowledgements}
This paper describes objective technical
results and analysis. Any subjective views or opinions that might be
expressed in the paper do not necessarily represent the views of the
U.S. Department of Energy or the United States Government. Sandia
National Laboratories is a multimission laboratory managed and
operated by National Technology and Engineering Solutions of Sandia
LLC, a wholly owned subsidiary of Honeywell International, Inc., for
the U.S. Department of Energy's National Nuclear Security
Administration under contract DE-NA-0003525. SAND2022-1366 O.

\section*{Appendix A}
This appendix provides details on the discretization of optimal control problems using adaptive pseudospectral methods. The fundamental ideas are common in fluid dynamics \cite{spectral_fluids_book} and became popular in trajectory planning thanks to the pioneering work of Fahroo and Ross \cite{fahroo_ross_ps_2002,psop_review_ross}. Adaptive pseudospectral methods have gained significant attention in recent years thanks to the work of Rao and collaborators \cite{hp_adap_psopt_rao,gpops}. 

For simplicity of the exposition we consider problems of the form 
\begin{align}\label{eq:OCP}
 \tag{OC} 
& \min_{\x,\u}  \int_0^{T} C_{run}(\x(t),\u(t),t)dt + C_{final}(\x(T),\u(T)) \\
& s.t. \nonumber \\
& \begin{dcases} \label{eq:dynamics_oc}
 \dot{\x}(t) = \f(t,\x,\u) \qquad t \in (0,T) \\
  \x(0) = \x_0
 \end{dcases}
\end{align}
and note that the subsequent developments may be easily extended for the more general case with inequality constraints, final time constrains, and a free final time. We have suppressed dependence on $\g$ for notational simplicity.

Pseudospectral methods discretize the dynamics~\eqref{eq:dynamics_oc} by representing the state and control in a finite dimensional basis and collocating the ODE system at a finite number of nodes in time. In particular, we consider a partition of the time interval $[0,T]$ into $p$ subintervals 
\begin{eqnarray} \label{eq:partition}
(0,T] = \bigcup_{i=1}^p (t_{i-1},t_i]
\end{eqnarray}
where $0=t_0 < t_1 < t_2 < \cdots < t_p = T$. We approximate the state variables using $N$ local (in the subintervals) Lagrange polynomials at Gauss-Lobatto nodes. Let $\{\mathcal Y_i\}_{i=1}^N$ denote this set of basis functions and $\y=(y_1^1,y_2^1,\dots,y_N^1,y_1^2,\dots,y_N^2,\dots,y_N^n)^T \in \r^{K_n}$, $K_n=nN$, denote the coordinate representation of the approximation, i.e.
\begin{eqnarray} \label{eq:state_approx}
x_k(t) \approx \sum\limits_{i=1}^N y_i^k \mathcal Y_i(t).
\end{eqnarray}

The controller is also discretized via an expansion in basis functions $\{\mathcal Z_j\}_{j=1}^M$; however, they may and in many cases are different than $\{\mathcal Y_i\}_{i=1}^N$. In particular, the state basis functions are adapted to resolve the dynamics whereas the control basis is designed to enforce the users desired smoothness in the control solution, although it may also be adapted to focus nodes in regions with fast time scales if the user desires. Let $\z=(z_1^1,z_2^1,\dots,z_N^1,z_1^2,\dots,z_N^2,\dots,z_N^n)^T \in \r^{K_m}$, $K_m=mM$, be the coordinates for the control, i.e.
\begin{eqnarray*}
u_k(t) \approx \sum\limits_{j=1}^N z_j^k \mathcal Z_j(t).
\end{eqnarray*}

We collocate the dynamics~\eqref{eq:dynamics_oc} by differentiating~\eqref{eq:state_approx} and evaluating the derivative at the time nodes. Enforcing that $\dot{\x}(t) = \f(t,\x,\u)$ at the time nodes yields the system of equations
\begin{eqnarray*}
\mathbf r_{coll}(\y,\z) = Y_t \y - \vec{\xi}(\y,\z) = \vec{0}
\end{eqnarray*}
where $\vec{\xi}(\y,\z)$ is a vector corresponding to stacking together evaluations of $\f(t,\x,\u)$ at the time nodes and $Y_t$ is a matrix populated with time derivatives of the basis functions $\{\mathcal Y_i\}$. The vector $\mathbf r_{coll}(\y,\z)$ contains $(N-p-1)n$ nonlinear equations as we do not collocate at $\{t_i\}_{i=0}^p$, the interface, initial, and terminal time nodes. Rather, the initial conditions are enforced directly yielding $n$ equations which we denote as $\vec{r}_{init}(\y,\z) \in \r^n$, and the terminal nodes on each subinterval $\{t_i\}_{i=1}^p$ are enforced by requiring that $\x(t_i)-\x(t_{i-1})=\int_{t_{i-1}}^{t_i} \f(t,\x,\u)$. This yields an additional $pn$ equations which we denote as $\vec{r}_{inter}(\y,\z) \in \r^{pn}$. Combining $\mathbf r_{coll}$, $\vec{r}_{init}$, and $\vec{r}_{inter}$ yields a system of $nN$ nonlinear equations
\begin{align*}
\mathbf r(\y,\z) = 
\left( \begin{array}{c}
\mathbf r_{coll}(\y,\z) \\
\vec{r}_{init}(\y,\z)  \\
\vec{r}_{inter}(\y,\z) \\
\end{array} \right) 
= \vec{0}
\end{align*}
whose solution is the coordinates for an approximate solution of the dynamics~\eqref{eq:dynamics_oc}. To ensure a quality approximation, the time interval partition~\eqref{eq:partition} is adapted to allocate time nodes in regions of faster dynamics.

The objective function is discretized by approximating the integral term
\begin{eqnarray*}
 \int_0^{T} C_{run}(\x(t),\u(t),t)dt \approx \sum\limits_{i=1}^N w_i c_{run}(\y,\z)
\end{eqnarray*}
where $\{w_i\}_{i=1}^N$ are the integration weights corresponding to applying Gauss-Lobatto integration on the subintervals and $c_{run}(\y,\z)$ is the evaluation of $C_{run}$ at the $i^{th}$ time node using the discretized state and control. We represent the final time objective $C_{final}(\x(T),\u(T))$ with $c_{final}(\y,\z)$ by evaluating $C_{final}$ using the coordinates for the final time node.

The discretized optimal control problem is
\begin{align} \label{eq:opt_control_discretized}
& \min_{\y,\z}   \sum\limits_{i=1}^N w_i c_{run}(\y,\z) + c_{final}(\y,\z)\\
& s.t. \nonumber \\
& \vec{r}(\y,\z) = 0 \nonumber 
\end{align}
and may be solved using standard nonlinear programming methods. 

\bibliographystyle{plain}
\bibliography{dasco}

\end{document}